\documentclass[english,11pt]{article}

\usepackage[a4paper, margin=2.1cm]{geometry}
\usepackage[english]{babel} 
\usepackage{mathtools}
\usepackage{times}
\usepackage{mathrsfs} 
\usepackage{amsmath}
\usepackage{amsthm}
\usepackage{amsfonts}
\usepackage{amssymb}
\usepackage[round]{natbib} 
\usepackage{thmtools}
\usepackage{relsize}
\usepackage{xcolor}
\usepackage[official]{eurosym}
\declaretheorem[style=definition,qed=$\triangle$]{example}

\newtheorem{theorem5}{Lemma}

\newtheorem{theorem1}{Definition}

\newtheorem{theorem4}{Corollary}

\newtheorem{theorem0}{Theorem}

\begin{document}

\title{Stable Spare Parts Pooling for Military Weapon Systems}
\author{Loe Schlicher, Marco Slikker, Willem van Jaarsveld\\
\footnotesize{School of Industrial Engineering, Eindhoven University of Technology} \\
\footnotesize{P.O. Box 513, 5600 MB, Eindhoven, The Netherlands}}

\date{\today}

\maketitle
\begin{abstract}
We study under which circumstances Departments of Defenses should be willing to deploy a joint parts part pooling program for their major weapon systems. 
Using cooperative game theory and Markov Decision Processes, we demonstrate that the type of pooling strategy plays a crucial role in the success of such a joint
spare parts pool. More precisely, we show that a joint spare parts pool may not last long –or even not arise–
if full pooling is applied, while it is stable under threshold pooling. 
\end{abstract}
\section{Introduction}

Various \emph{joint spare parts pooling programs} for common major weapon systems, such as the NH90 helicopter, F-16 and F-35 combat aircraft, have been initiated in the last decades. These 
pooling programs, however, have proven to be difficult to implement. As discussed in the RAND report of \cite{lorell2016review}  "[...] one of the most significant barriers to spare parts pooling [...] is the question of \emph{security of supply}." Put simply, that is "[...] the level of assurance that a specific spare part will be received when it is needed". 
Another significant barrier to spare parts pooling is "[...] the risk
of shirking—the situation in which a partner is unwilling to pay its allocated share of the total cost. 
As an example, the RAND report discusses the case of the USAF, who identified shirking as a major potential risk for participation in the F-35 spare parts pooling initiative. 

Given the \emph{advantages} of spare parts pooling, such as (i) the reduction in investment and storage costs as well as (ii) the improved availability of spare parts and so increased readiness of military weapon systems (see, e.g., \citet{EDA})
, it is of utmost importance to think about ways to mitigate the aforementioned two barriers. 
For the first barrier, the security of supply, it seems to be important to \emph{inform} partners about the actual likelihood of out-of-stock situations in the pool and the fact that this probability depends on an upfront \emph{collective} decision on the \emph{number} of spare parts to stock and the \emph{type of pooling strategy} (i.e., who is allowed to demand when and what?). For the second barrier, the risk of shirking, it seems to be important to 
 \emph{allocate} the total costs of the joint spare parts pool in such a way that no partner has reasons to leave the joint spare parts pool. Leaving a spare parts pool seems to be a serious concern; the RAND report of \cite{lorell2016review} mentions that for the F-35 program "[...] language was included to allow participants to opt out of the joint spare parts pool, should a nation desire to establish (and be willing to pay for) a separate spare stock".

In this paper, we will study the \emph{collective decision} of how many \emph{spare parts to stock}, what \emph{type of pooling strategies} to use and subsequently how to \emph{allocate the costs} 
when collaborating Departments of Defenses (DoDs) are willing to deploy a joint parts pooling program. 
We consider two types of pooling strategies in the paper, namely \emph{full-pooling} and \emph{threshold pooling}, which should be interpreted as a set of rules describing 
 when DoDs are allowed to demand for spare parts. Full-pooling, for instance, prescribes that \emph{everyone} can demand for parts as long as some spare parts are available in the pool. Threshold pooling, on the other hand, prescribes that DoDs can \emph{only} demand for spare parts if the number of spare parts left in the spare pool is above a certain threshold level. This threshold could, for instance, be based on unique characteristics of the involved DoDs, such as (i) the internal costs for \emph{holding} spare parts and (ii) the direct and indirect (cascade effect) costs associated to equipment \emph{downtime}, as a result of receiving no spare parts. 
  We use \emph{Markov Decision Processes} (MDP), a well-known method in the spare parts literature (see, e.g. \cite{wong2005stocking,wong2007cost, kranenburg2009new}), to model the two types of pooling strategies and
  to study the collective decision of how many spare parts to stock for both types of pooling strategies, using a long-term average cost criterion.

We address the subsequent question of \emph{cost allocation} using \emph{cooperative game theory}. More precisely, 
we introduce two \emph{cooperatives games}, one for each type of  pooling strategy. In both games, the same group of DoDs form a joint spare parts pool and we refer to this group of DoDs as the \emph{grand coalition}. Using MDPs, we determine --for full and threshold pooling, respectively--  the optimal number of spare parts to stock 
 and their associated long-term average costs. 
 On top of this, we also calculate the long-term average cost for all possible subgroups of DoDs, which we refer to as \emph{coalitions}. 
 Subsequently, 
 we investigate whether the costs of the grand coalition can be divided amongst the DoDs in such a way that no coalition has financial incentives to split off from the grand coalition. This concept of (groups of) DoDs having the possibility to split off seems to be realistic (see, e.g., the RAND report of \citet{lorell2016review} to which we referred earlier).
 
 In cooperative game theory, the set of payoff vectors allocating the costs of the grand coalition in such a way that no coalition, consisting of a  single or multiple DoDs, has financial incentives to split off is called the \emph{core} (see, e.g., \cite{peleg2007introduction}). Consequently, for our two cooperative games, we are basically
 interested in the question whether their cores are non-empty, i.e., whether the core is non-empty under full and threshold pooling. The outcomes turn out to be quite delicate:
we prove that if DoDs 
 apply \emph{threshold pooling}, the core is always non-empty. However, if  DoDs decide to apply \emph{full pooling}, core non-emptiness is not guaranteed any longer. More precisely, core non-emptiness is violated in situations where heterogeneous demand and cost structures arise. In such situations, DoDs with relative low downtime costs will typically overdemand from the spare parts pool, while nobody is willing too pay for that in the end. With threshold pooling, these DoDs will be limited to overdemand from the joint spare parts pool. These insights are practically relevant, because current spare parts pooling initiatives rely on full pooling (see, e.g., \cite{lorell2016review}).


We would like to emphasize that 
the coalitional costs in our threshold pooling game arise as the costs of an MDP, and as a consequence there is no closed-form nor implicit form available that describes our coalitional values in general. Hence, it is unlikely to prove core non-emptiness in terms of the primal characteristics of our underlying spare parts pooling situation. A common alternative to prove core non-emptiness is showing that the game is \emph{concave}, but unfortunately our game is not concave in general. Instead, we demonstrate that threshold pooling of a balanced collection of coalitions can be \emph{merged} to a feasible threshold pooling strategy for the grand coalition. In doing so, we use that our underlying situation can be described by an MDP and threshold pooling as
a stationary policy in this MDP. This insight can be leveraged to establish the \emph{balancedness conditions} of \citet{bondareva1963some} and \citet{shapley1967balanced}, which enables a proof that the game is balanced. In turn, this implies that the core of our threshold pooling game is non-empty (see, e.g., \citet{bondareva1963some} and \citet{shapley1967balanced}). 
Using MDPs in order to prove balancedness is a relatively new technique that was first proposed in \cite{schlicherMOR}. 
This technique is also applied to our game, that involves a one-time decision about the number of spare parts to stock, next to identifying optimal threshold levels.
This application, however, brings two new challenges, the first of which is how to deal with the one-time
decision on how many spare parts to stock.  
Second, the technique presented in \cite{schlicherMOR} is applicable to MDP formulations with a \emph{one-dimensional}
state space only, while heterogeneous holding costs and storage capacities --as present in our setting-- call for a \emph{multi-dimensional} state space approach.

Our work contributes to the stream of literature using cooperative game theory to analyze spare parts pooling situations. Pioneer work in this field is done by \cite{wong2007cost}. They consider a setting where various parties can collaborate by sharing spare parts via lateral transshipments. Sharing is only allowed if a parties' own spare parts stock is above a certain threshold, which could be seen as a variant of threshold pooling. The authors introduce four allocation rules and investigate numerically whether their payoff vectors belong to the core for \emph{one} specific setting only. \cite{akkerman20222} demonstrates, however, that core non-emptiness cannot be guaranteed in general. This also holds for another spare parts pooling game proposed by \cite{guajardo2015cost}, where parties pool spare parts via setting service level targets jointly.


Next to these numerically oriented papers, there also exist some studies that focus on proving core non-emptiness for specific classes of spare parts pooling situations. An example in that regard is \cite{karsten2012inventory} who analyzes a similar spare parts pooling game as \cite{wong2007cost}. However, they assume that spare parts are always shared amongst participants (i.e., they apply full-pooling). Under rather restrictive assumptions, such as homogeneous demand rates, downtime costs and holding costs, \cite{karsten2012inventory} prove that the core is always non-empty. However, in the absence of some of these assumptions (e.g., a setting with heterogeneous downtime costs), the core can be empty. For another, rather limited, class of spare parts pooling situations, namely those where each player stocks exactly one item, but with heterogeneous downtime costs, \citet{schlicher2017probabilistic} demonstrates that core non-emptiness is guaranteed. Our paper contributes to the theoretical stream, by also proving core non-emptiness for a class of joint spare parts pooling situations. 

Cooperative game theory has been successfully applied to various types of real-life collaborative settings, including inventory pooling (\cite{meca2004inventory, guardiola2009production, fiestras2011cooperative}), transportation (\cite{agarwal2010network, algaba2019horizontal, van2021joint, schlicher2022stable}), supply chains (\cite{hall2010capacity, zhang2009cost, leng2009allocation}), retail (\cite{anily2014subadditive, chen2016duality}), health care (\cite{westerink2020core}) and the (high-tech) service industry (\cite{anily2010cooperation, karsten2015resource}). Our paper contributes to the  'inventory pooling' and 'service industry' settings.

The remainder of this paper is organized as follows. In 
Section~\ref{sec:preliminaries}, we discuss preliminaries
on cooperative
game theory as well as on MDPs. Subsequently, we introduce the joint spare parts pooling situation and associated games in Section~\ref{sec:model}. Our main result that threshold pooling game have a non-empty core is presented in Section~\ref{sec:mainresult}. We conclude in Section~\ref{sec:conclusion}. All proofs of lemmas are relegated
to the appendix.

\section{Preliminaries}\label{sec:preliminaries}

In this section, we provide some basic elements of discrete time Markov decision processes as well as cooperative game theory. Readers familiar with these topics can move to Section \ref{sec:model} directly.


\subsection{Discrete Time Markov Decision Processes}
\label{section:MDP}
In this section we present some basic concepts of discrete time Markov decision processes (MDPs).

To model sequential decision under uncertainty, one can make use of the mathematical framework of MDPs.
Let us consider a set $T = \mathbb{N} \cup \{0\}$ of time epochs, a countable set $\mathscr{Y}$ of states, a finite set $\mathscr{A}(y)$ of actions for each $y \in \mathscr{Y}$, non-negative costs $\mathscr{C}(y,a)$ for each $y \in \mathscr{Y}$ and all $a \in \mathscr{A}(y)$, and transition probabilities $p(y' \vert y,a)$ for all $y' \in \mathscr{Y}$, all $y \in \mathscr{Y}$, and all $a \in \mathscr{A}(y)$ with $\sum_{y' \in \mathscr{Y}} p(y' \vert y,a) = 1$ for all $y \in \mathscr{Y}$ and all $a \in \mathscr{A}(y)$. Then, tuple $(T,\mathscr{Y},\mathscr{A},\mathscr{C},p)$ with $\mathscr{A} = (\mathscr{A}(y))_{y \in \mathscr{Y}}$, $\mathscr{C} = (\mathscr{C}(y,a))_{y \in \mathscr{Y}, a \in \mathscr{A}(y)}$ and $p = (p(y' \vert y,a))_{y',y \in \mathscr{Y}, a\in \mathscr{A}(y)}$ is called a discrete time Markov decision process. 

Now, let $t \in T$ be a time epoch. A \emph{decision rule} $\omega_t = (\omega_t(y))_{y \in \mathscr{Y}}$ indicates for all states $y \in \mathscr{Y}$ which action to choose at time epoch $t$. Next, a \emph{policy} $\omega = (\omega_t)_{t \in T}$ is a sequence of decision rules for all time epochs. Let $X_t$ with $t \in \mathbb{N}_0$ be a random variable indicating the state at time $t$. Note, $X_t$ depends on $\omega$ and $X_0$. If initially $X_0 =y \in \mathscr{Y}$, the \emph{long-run average costs per time epoch} under policy $\omega$ are
\begin{equation*}J_{\omega}(y) = \lim \sup_{n \to \infty} \frac{1}{n} \mathbb{E}_{\omega} \left[ \sum_{t=0}^{n-1} \mathscr{C}(X_t,\omega_t(X_t)) \vert X_0=y\right].\end{equation*}

We define $\Omega$ to be the set of all policies and let $J^*(y) = \inf_{\omega \in \Omega} J_{\omega}(y)$ for all $y \in \mathscr{Y}$. There exists a class of MDPs for which there exists a constant $J^*$ such that $J^* = J^*(y)$ for all $y \in \mathscr{Y}$. For this specific setting, $J^*$ is defined as the \emph{minimal long-run average costs per time epoch} and policy $\omega \in \Omega$ is called \emph{optimal} if $J_{\omega}(y) = J^*$ for all $y \in \mathscr{Y}$. Moreover, a policy is called \emph{stationary} if there exists an $f$ such that $\omega_t = f$ for all $t \in T$. In addition, we denote such a stationary policy by $f = (f(y))_{y \in \mathscr{Y}}$.

For an MDP, the \emph{value function} $V_t$ is defined by
\begin{equation}V_{t+1}(y) = \min_{a \in \mathscr{A}(y)} \left\{ \mathscr{C}(y,a) + \sum_{y' \in \mathscr{Y}} p(y' \vert y,a) \cdot V_t(y')\right\} \mbox{ for all } t \in T \mbox{ and all } y \in \mathscr{Y},\label{valuefunctiondef} \end{equation}
\noindent with $V_0(y) = 0$ for all $y \in \mathscr{Y}$.


%


\noindent The following result states that under two conditions, the minimal long-run average costs per time epoch exist, are attained under a stationary policy and moreover, coincide with the limit of the value function divided by the number of time epochs, when time goes to infinity. The first condition refers to   \emph{irreducible} Markov chains and \emph{positive recurrent} Markov chains. A Markov chain, which is attained under a given policy, is said to be \emph{irreducible} if for all states $y \in \mathscr{Y}$ it holds that it can be reached from each other state $y' \in \mathscr{Y} \backslash \{y\}$. A Markov chain is said to be \emph{positive recurrent} if for all states $y \in \mathscr{Y}$ it holds that the expected return time (to state $y$) is finite. The second conditions relates to the fact that the state space needs to be finite. We now present the theorem, which is derived in \citet[proposition 4.3]{sennott1996convergence}.



\begin{theorem0} \label{sennot} Let $(T,\mathscr{Y},\mathscr{A},\mathscr{C},p)$ be an MDP. If $(i)$ there exists a stationary policy $f$ inducing an irreducible and positive recurrent Markov chain on $\mathscr{Y}$, and satisfying $J_{f}(y) < \infty$ for all $y \in \mathscr{Y}$, and $(ii)$ there exists an $\epsilon >0$ such that $\{y \in \mathscr{Y} \vert \exists a \in \mathscr{A}(y) : \mathscr{C}(y,a) < J_{f}+ \epsilon\}$ is finite, then
\begin{equation} J^* = \lim_{t \to \infty} \frac{V_t(y)}{t} \mbox{ for all } y \in \mathscr{Y},  \end{equation}

\noindent and moreover, there exists an optimal stationary policy. \end{theorem0}

In this paper, we restrict our attention to MDPs with finite state spaces. This implies that it becomes superfluous to check whether an irreducible Markov chain is positive recurrent as every irreducible Markov chain is positive recurrent in a finite state space
(see, e.g., \cite[Theorem 5.71 (ii)]{modica2012first}. Moreover, the second condition of Theorem \ref{sennot} is always satisfied when $\mathscr{Y}$ is finite.

\subsection{Cooperative Game Theory}\label{sec:cooperativegametheory}

A \emph{cooperative cost game with transferable utility}, shortly called \emph{game}, can be described by a pair $(N,c)$ with $N$ a finite set of \emph{players} $N = \{1,2,\ldots,n\}$ and a \emph{characteristic function} $c: 2^N  \to \mathbb{R}$. 
A \emph{coalition} is a subset $S \subseteq N$ with $c(S)$ representing the costs incurred by the players in $S$. It is assumed that the costs can be transferred freely among the players in $S$. Moreover, the set $N$ is called the \emph{grand coalition}. A \emph{cost vector} for a game $(N,c)$ is a vector $x \in \mathbb{R}^N$ describing how to allocate the costs, where player $i \in N$ is allocated $x_i$. On top of this, a cost vector $x \in \mathbb{R}^N$ is called \emph{efficient} if $\sum_{i \in N}x_i = c(N)$, i.e., all costs are distributed amongst the players of the grand coalition $N$. A cost vector $x \in \mathbb{R}^N$ is called \emph{stable} if no group of players has an incentive to leave the grand coalition $N$ (i.e., $\sum_{i \in S}x_i \leq c(S)$ for all $S \subseteq N$). The set of efficient and stable cost vectors of $( N,c)$ is called the \emph{core} of $(N,c)$ and denoted by $core(N,c)$. 

\subsection{Balancedness of games}\label{sec:balancedness}
In this paper, we make a link between core non-emptiness and \emph{balancedness} of cooperative games, which we summarize next. A map $\kappa : 2^N \backslash \{\emptyset\} \to [0, 1]$ is called a \emph{balanced map} for $N$ if
$$ \sum_{S \in 2^N : \hspace{1mm} i \in S} \kappa_S = 1 \hspace{2mm} \mbox{ for all } i \in N.$$

A collection $\mathscr{B} \subseteq 2^N \backslash \{\emptyset\}$ is called \emph{balanced} if there exists a balanced map $\kappa$ for which $\kappa_S > 0$ for all
$S \in \mathscr{B}$ and $\kappa_S=0$ otherwise. In addition, a collection $\mathscr{B} \subseteq 2^N \backslash \{\emptyset\}$ is called \emph{minimal} \emph{balanced} if there exists no proper subcollection of $\mathscr{B}$ that is balanced as well. As discussed in \citet{peleg2007introduction}, for every minimal balanced collection $\mathscr{B} \subseteq 2^N \backslash \{\emptyset\}$ there exists exactly one associated balanced map $\kappa$. Moreover, for such a balanced map it holds that $\kappa_S \in \mathbb{Q}$ for all $S \in \mathscr{B}$ (\citet{norde2002dual}). 

A game $(N, c)$ is called \emph{balanced} if for every minimal
balanced collection $\mathscr{B} \subseteq 2^N \backslash \{\emptyset\}$ with associated balanced map $\kappa$ it holds that
$$\sum_{S \in \mathscr{B}} \kappa_S \cdot c(S) \geq c(N).$$

Independently from each other, \citet{bondareva1963some} and \citet{shapley1967balanced} provided a sufficient and necessary conditions for core non-emptiness, which we present next.


\begin{theorem0} A game $(N, c)$ is balanced if and only if its core is non-empty. \label{thm:degrotestelling} \end{theorem0}

Let $\mathscr{B} \subseteq 2^N \backslash \{\emptyset\}$ be a minimal balanced collection and $\alpha \in \mathbb{N}$ be the smallest integer for which $\kappa_S \cdot   \alpha \in \mathbb{N}$ for all $S \in \mathscr{B}$. Note that such an $\alpha$ exists as $\kappa_S \in \mathbb{Q}$ for all $S \subseteq N$. As a shorthand notation, we use $b_S = \alpha \cdot \kappa_S$ for all $S \in \mathscr{B}$, suppressing the dependency of $\alpha,b_S$, and $\kappa_S$ on $\mathscr{B}$. Now, by applying the result of Theorem \ref{thm:degrotestelling}, it suffices to prove balancedness by verifying the inequality 
\begin{equation} \label{balanced} \sum_{S \in \mathscr{B}} b_S \cdot c(S) \geq \alpha \cdot c(N), \end{equation}
 for all minimal balanced collections $\mathscr{B} \subseteq 2^N \backslash \{\emptyset\}$.
\section{Joint spare parts pooling situation and two associated games}\label{sec:model}

In this section, we define the joint spare parts pooling situation and two associated
games. Besides, we describe how the coalitional values of our games relate to the value function of a corresponding MDP. While introducing the joint spare parts pooling situation, we use the F-35 combat aircraft as a motivational example. 



\subsection{Joint spare parts pooling situation}
We consider an environment with a set of  $N \subseteq \mathbb{N}$ Departments of Defence (DoDs),  each operating a large fleet of the same type of military asset (e.g., the F-35 aircraft). 
We focus on \emph{one specific repairable part} of this asset only (e.g., the engine of a F-35) that may be critical to the asset's functioning.
We assume that the specific part is subject to failure according to a Poisson distribution with rate $\lambda_i \in \mathbb{R}_{>0}$ for all $i \in N$, which is common in literature on  military spare parts operations \citep{basten2014system}. The assumption of a constant failure rate is reasonable if the fleet size is sufficiently large. This is, for instance, the case for the F-35 combat aircraft with more than 2,400 F-35 in the US \citep{behoefteNL} and 550 in Europe \citep{newslockheed}. Next, we assume that each DoD has a \emph{local stock point} close to its  asset fleet with fixed storage capacity  $C_i \in \mathbb{N}_{\geq 0}$. 
As turns out later, this assumption is not restrictive. We assume that the \emph{cost for holding} a single spare part at local stock point of DoD $i \in N$ is denoted by $h_i \in \mathbb{R}_+$ per time unit. That is, we allow for cost heterogeneity because insurance and capital costs may differ widely
between DoDs. For ease of exposition we assume that DoDs are ordered such that $i>j$ implies that $h_i\ge h_j$, i.e. DoDs with a lower index have lower holding costs. Finally, we assume that the local stock points of the DoDs are geographically close to each other, such that the time to transport parts between the points is within hours. 
For the F-35 combat aircraft setting this would, for instance, be the case if some European countries such as the Netherlands, Germany, Belgium and Denmark would team up and form a joint spare pool. Another example could be a setting where various DoDs together form a large base camp in a mission area for several years (see, e.g., \cite{Defensie}).

We assume that DoDs decide to collaborate by forming a \emph{joint spare parts pool}. 
This pool is virtual in the sense that each DoD still stocks spare parts \emph{locally} (i.e., stocks spare parts at their local stock point), but the allocation of the spare parts is assumed to be \emph{centrally} organized. That is, there is an authority who dictates which DoDs are allowed to take parts from the joint spare parts pool whenever there is demand. Before the DoDs operationalize the joint spare parts pool, they
collectively decide upon the \emph{total} number of spare parts to stock and where to locate them initially, while taking into account 
the local storage capacities per DoD.


 Once the pool has been operationalized, the authority has to accept or reject demand when it arrives.
 If the authority allows the DoD to demand, it will be taken from the local stock point --if available-- and otherwise from one of the neighboring local stock points. 
 Next, the DoD sends the failed, but repairable, part to the Original Equipment Manufacturer (OEM), who is certified to  repair it. The OEM repairs defective parts one-by-one and we assume the \emph{repair times} to be independent and identically distributed according to an exponential distribution with rate $\sum_{i \in N}\mu_i$, where $\mu_i \in \mathbb{R}_+$ for all $i \in N$. Rate $\mu_i$ could, for instance, be interpreted as the repair capacity reserved for DoD $i \in N$ in a long-term contract negotiated with the OEM. Once the repair completes, the part is returned to the joint spare parts pool, i.e., to the local stock point with lowest holding cost having capacity available.


If the authority does not allow the DoD to demand, or when there are no spare parts left in the pool, the DoD will instigate an emergency procedure: a 'new' spare
part is obtained from the responsible OEM to replace the failed part directly. We assume that the time to transport this new spare part is significantly longer than the time to transport a spare part between local stock points. This assumption seems to be reasonable for the F-35 case, where Lockheed Martin could be the responsible OEM to start up an emergency procedure. In that case, most likely a spare part will be shipped from the United States to Europa, which typically takes much longer than transferring spare parts between local stock points in Europe. Finally, we assume that the failed part is collected by the OEM, but not returned later on. As a consequence, the total
number of spare parts that is in repair at the OEM (for all DoDs) or in the local stock points,
remains constant over time. In Table \ref{table1}, we summarize the two decisions the authority can make as  well as the associated consequences.

\begin{table}[h!]
    \centering
    \begin{tabular}{ c | c c c c} \hline
    decision  &   spare part taken from &  failed part sent to & repaired by & repaired part\\ \hline
     \emph{accept}  &   spare parts pool &  OEM & OEM &  goes back to spare parts pool \\
    \emph{reject} &    OEM (emergency) & OEM &  - & stays at OEM \\ \hline
    \end{tabular}
    \caption{Two possible decisions once a part has failed}
    \label{table1}
\end{table}

If a spare part is not served from the joint spare parts pool (see option 1 of Table \ref{table1}), a DoD has to wait for a spare parts via the emergency procedure, implying that the underlying military asset cannot function properly for a significant amount of time. To account for this equipment downtime, we introduce a \emph{downtime cost} $d_i \in \mathbb{R}_{\geq 0}$ per DoD $i \in N$, reflecting the direct and indirect losses arising from the inability of the military asset to respond to operational demands. Over the years, various DoDs have developed quantitative tools to assess these direct and indirect costs of downtime of military assets (see, e.g., \cite{fuerst1991model} \cite{peltz2002diagnosing} and \citet{nato2003cost}). 
Next to the direct and indirect costs of equipment downtime, $d_i$ could also reflect the cost to facilitate an emergency procedure, such as contractual costs associated to an emergency procedure set by the OEM. 
Please, observe that we do not account for downtime costs if a spare part is served from the joint spare parts pool. We neglect this cost, because in practice it might already take hours to replace a disfunctioning part from an asset, which is more or less the same time to transport a spare part between DoD stock locations such as the Netherlands and Germany or Belgium. Hence, in practice there are costs associated to downtime, but these are unavoidable and so neglected in our model.

Finally, we assume the DoDs' objective is to minimize the long-run average costs. They do so by collectively deciding upon the total number of spare parts to stock and where to locate them initially, while taking into account the local storage capacities per DoD. On top of that, via the authority, they consider two 'accept or reject' strategies. In the first one, which we refer to as \emph{full-pooling}, the authority accepts \emph{all} demand. In the second one, the authority accepts or rejects demand such that the long-run average costs is minimized. According to 
\cite{ha2000stock} 
, an optimal 'accept or reject 'structure has the form of a threshold strategy.\footnote{The stock rationing problem of a single-item, make-to-stock production system with several demand classes and lost sales studied by \cite{ha2000stock} is mathematically equivalent to our spare parts pooling model when production time is exponentially distributed. 
In Theorem 1 of \cite{ha2000stock} the author shows optimality of a threshold policy under a discounted cost criterion and argues that this result can be extended to the long-run average cost criterion (p.80, l.-1 to l.-7).} That is, the authority should only accept demand of a certain DoD if the total number of spare parts left in the pool is above a certain DoD specific threshold. In this way, some spare parts are reserved for the DoDs with relatively high downtime costs.\footnote{The policy-iteration algorithm (see, e.g., \cite[p.247]{tijms2003first}) can, for example, be employed to actually find optimal thresholds levels and the optimal number of spare parts to stock when the joint spare parts pooling situation is modelled as an MDP.} In the paper, we refer to such pooling strategies as \emph{threshold pooling}.


To analyse this setting, we define a \emph{joint spare parts pooling situation} as a tuple $(N,\lambda, C, h, \mu, d)$ with $N$, $\lambda = (\lambda_i)_{i\in N}$, $C=(C_i)_{i \in N}$, $h= (h_i)_{i \in N}$, $\mu = (\mu_i)_{i \in N}$, and $d=(d_i)_{i \in N}$ as defined above. For short, we use $\theta$ to refer to a joint spare parts pooling situation and specify the elements of $\theta$ by $(N, \lambda, C, h,\mu,d)$ unless stated differently. Moreover, we use $\Theta$ to refer to the set of joint spare parts pooling situations. 


\subsection{Full Pooling and Threshold Pooling Games}

For every joint spare parts pooling situation $ \theta \in \Theta$, we introduce and analyze two associated cooperative games, namely the full-pooling game $(N, c^{\theta}_F)$ and the threshold-pooling game $(N,c^{\theta}_T)$. In these two games, the characteristic cost value $c^{\theta}_F(S)$ and $c^{\theta}_T(S)$, respectively, reflect the long-run average costs of the joint spare parts pooling situation restricted to the players (i.e., DoDs) in coalition $S \subseteq N$. That is, we only consider $(\lambda_i)_{i \in S}$, $(C_i)_{i \in S}$, $(h_i)_{i \in S}$, $(\mu_i)_{i \in S}$ and $(d_i)_{i \in S}$. In doing so, we assume that the repair capacity reserved for the DoDs in a long-term contract negotiated with the OEM reads $\sum_{i \in S} \mu_i$ in coalition $S \subseteq N$. 
This implies that smaller coalitions thus have smaller bargaining power while negotiating contractual agreements.

\newcommand{\bs}{\xi}

\subsection{MDP Formulation}\label{sec:MDP}

 The coalitional values for both the full-pooling game and the threshold pooling game can be described by a (discrete time) MDP. This is allowed since the decision problem per coalition can be recognized as a semi-Markov decision problem, which can be converted to an equivalent MDP by applying \emph{uniformization} (see, e.g., \citet{lippman1975applying}). For that, we add fictitious transitions of a state to itself to ensure that the total rate out of a state is equal for all states, the so-called \emph{uniformization
rate}. Then, we consider the embedded discrete-time MDP by looking at the system only at transition instants, which occur according to a Poisson process, with
as rate the uniformization rate. 
In what follows, we first describe the MDP for the threshold pooling game. Let $\theta \in \Theta$, $S \subseteq N$ with $S \not = \emptyset$ and assume that DoDs collaborate via threshold pooling.

\subsubsection{State and action space of threshold pooling game}\label{sec:statespace}

\noindent 

We define the state space to be $\mathscr{Z}^S = \prod_{i \in S}\{0,1,\ldots,C_i\}$ with state $z=(z_i)_{i \in S} \in \mathscr{Z}^S$ describing the the number of spare parts in stock per local stock point, i.e., $z_i$ is the number of spare parts in local stock point $i$ for all $i \in S$.
Moreover, for all $z \in \mathscr{Z}^S$ we introduce action sets \begin{equation*}
\begin{aligned}
\mathscr{A}^{S}_{-}(z) &= \left\{(x_i)_{i \in S} \hspace{1mm} \Bigg\vert \hspace{1mm}\begin{matrix} x_i \in \{0,\min\{1,z_i\}\} \hspace{1mm} \mbox{ for all } i \in S  \\ \sum_{i \in S} x_i \leq 1 \end{matrix}  \right\}; \\
\mathscr{A}^S_{+}(z) &= \left\{(x_i)_{i \in S} \hspace{1mm} \Bigg\vert \hspace{1mm}\begin{matrix} x_i \in \{0,\min\{1,C_i-z_i\}\} \hspace{1mm} \mbox{ for all } i \in S  \\ \sum_{i \in S} x_i \leq 1 \end{matrix}  \right\}. 
\end{aligned}
\end{equation*}
The actions in the first action space should be recognized as the possibilities for players to demand for spare parts, based on the current spare parts available, while the actions in the second action space should be recognized as the possibilities to sent back repaired spare parts to a local spare stock. 

For each state $z \in \mathscr{Z}^S$, we define action $a = ((a_{ij}^{-})_{j \in S}, (a_{ij}^{+})_{j \in S})_{i \in S}$. Here, $(a_{ij}^{-})_{j \in S} \in \mathscr{A}^S_{-}(z)$ and $(a_{ij}^{+})_{j \in S} \in \mathscr{A}^S_{+}(z)$ for all $i\in S$. The action space becomes $\mathscr{A}^S(z) = \prod_{i \in S} \mathscr{A}^S_{-}(z) \times  \prod_{i \in S} \mathscr{A}_{+}^S(z)$.

\subsubsection{Uniformization and state transitions of threshold pooling game}\label{sec:statetransitions}
Let $\gamma = \sum_{i \in N} \left[\lambda_i + \mu_i\right]$.  We will use $\gamma$ as the uniformization rate, which is independent
of $S$.   
Let $\lambda^*_i = \lambda_i / \gamma$, $\mu^*_i = \mu_i / \gamma$ and $h^*_i = h_i / \gamma$ for all $i \in N$. Now, $C^S(z,a)$ denotes the expected costs collected over a single (uniformized) time epoch, given that the system begins the period in state $z \in \mathscr{Z}^S$ and action $a = ((a_{ij}^-)_{j \in S},(a_{ij}^+)_{j \in S})_{i \in S} \in \mathscr{A}^S(z)$ is taken. For our situation, we have
\begin{equation*} C^S(z,a) =  \sum_{i \in S} \lambda^*_i \cdot \left(1-\sum_{j \in S}a_{ij}^-\right) \cdot d_i + \sum_{i \in S }h^*_i \cdot z_i\hspace{7mm} \mbox{ for all } z \in \mathscr{Z}^S \mbox{ and all } a \in \mathscr{A}^S(z).\end{equation*}

\noindent In addition, let $p^S(z' \vert z,a)$ denote the one-stage transition probability from state $z \in \mathscr{Z}^S$ to $z' \in \mathscr{Z}^S$ under action $a \in \mathscr{A}^S(z)$:
\begin{equation*} p^S(z' \vert z,a) = \left\{\begin{matrix}   \lambda^*_i \cdot a_{ij}^-&& \mbox{ if } && z' = z-e^S(j) \mbox{ for all } i \in S; \\
 \mu^*_i \cdot a_{ij}^+  && \mbox{ if } && z' = z+e^S(j) \mbox{ for all } i \in S;  \\
   1 - \mathlarger{\sum}_{i \in S} \left[\lambda^*_i \cdot \sum_{j \in S }a_{ij}^- + \mu_i^* \cdot \sum_{j \in S }a_{ij}^+   \right]  && \mbox{ if } && z' =z; \\
0 && \mbox { else,} \end{matrix} \right.\end{equation*}

\noindent for all $z \in \mathscr{Z}^S$ and all $a \in \mathscr{A}^S(z)$, with
$e^S(i) = (e^S_j(i))_{j \in S}$ for all $S \subseteq N$ and all $i \in S$, with $e^S_j(i) = 1$ if $i=j$ and $e^S_i(j)=0$ otherwise. 

We denote the value function in our embedded MDP for coalition $S$ at time $t+1$ by $W_{t+1}^S$. By some elementary rewriting of the value function of (\ref{valuefunctiondef}), we end up with the formulation below. 

\begin{theorem5} \label{grotemdp} Let $\theta \in \Theta$ and $S \subseteq N$. Then, for all $t \in \mathbb{N}_0$ and all $z \in \mathscr{Z}^S$, the value function satisfies
\begin{equation*} \label{lemma:lemma1}\begin{aligned} W_{t+1}^S(z) =& \sum_{i \in S} \left[\lambda_i^* \min_{x \in \mathscr{A}^S_{-}(z)} \left\{ W_t(z-x) + (1 - \vert\vert x \vert\vert_{_{1}})d_i\right\} + \mu_i^* \min_{x \in \mathscr{A}_{+}^S(z)} \left\{W_t(z+x) \right\} \right]  \\
&+ \left(1- \sum_{i \in S} \left[ \lambda^*_i + \mu^*_i\right] \right) \cdot W^S_t(z) + \sum_{i \in S} h_i^* \cdot z_i\end{aligned} \end{equation*}with $W^S_0(z)=0$ for all $z \in \mathscr{Z}^S$.
\end{theorem5}

The formulation of the value
function in Lemma \ref{lemma:lemma1} can be interpreted in the following way. With probability $\lambda_i^*$
, there is a demand arrival. Subsequently, the authority has to decide whether they want to accept this demand --and if so from which stock location-- or reject it. If the authority decides to reject it (i.e., $\vert \vert x \vert \vert_1 =0$)
a downtime cost $d_i$ is incurred and a transition back to the same state. If the authority decides to accept demand (i.e., $\vert \vert x \vert \vert_1 =1$), there is no cost involved, but there is a transition to
a state where exactly one of the local stock points has one spare part less. With probability $\mu_i^*$ there is a 
repair completion.
Subsequently, the authority has to decide whether they want to accept this completion --and if so where it goes--  or reject it. 
In both cases, there are no costs involved, but if the authority decides to reject it  (i.e., $\vert \vert x  \vert \vert_1 = 0)$, there will be a transition to the current state, while an acceptance (i.e., $\vert \vert x  \vert \vert_1 = 1)$ leads to a transition to a state where exactly one local stock point has one spare part more. Note, there is no possibility to accept a completion in state  $\vert \vert z \vert \vert_1 = C_i$, resembling the setting where all (repaired) spare parts are in stock. With probability
$1 - \sum_{i \in S}(\lambda_i + \mu_i) \geq 0$, 
there is a dummy transition back to the current state, ensuring that the probabilities sum to 1.

Finally, we define $g^S$ as the minimal long-run average costs \emph{per time epoch} of the MDP. There is a direct relation between $g^S$ and the original minimal long-run average costs per time unit of coalition $S$.


\begin{theorem5} 
\label{cSlimrelation}Let $\theta \in \Theta$ and $S \subseteq N$ with $S \not = \emptyset$. Then,
$$ c^{\theta}(S) = \gamma \cdot g^S = \gamma \cdot \lim_{t \to \infty} \frac{W_t^S(z)}{t} \mbox{ for all } z \in \mathscr{Z}^S.$$ 
\end{theorem5}

\subsubsection{MDP for full-pooling}
The MDP formulation for full-pooling is equivalent to the one for threshold pooling except that this time demand is always accepted (if inventory is available), i.e.,  action set $\mathscr{A}_{-}^S(z)$ for any $z \in \mathscr{Z}^S$ 
 is updated to
 \begin{equation*}
\begin{aligned}
 \left\{(x_i)_{i \in S} \hspace{1mm} \Bigg\vert \hspace{1mm}\begin{matrix} x_i \in \min\{1,z_i\} \hspace{1mm} \mbox{ for all } i \in S  \\ \sum_{i \in S} x_i \leq 1 \end{matrix}  \right\}. \end{aligned}
\end{equation*}


\subsubsection{Two examples: the full-pooling game and threshold pooling game}
In this final subsection, we discuss two examples. In the first example, we present a full-pooling game, while in the second game, we present a threshold pooling game. For both of them, we study the core. 
Given that we consider relatively small examples, it is possible to verify all threshold strategies by hand, but one can also use the policy-iteration algorithm to find the optimal number of spare parts to stock and threshold levels.

\begin{example}
\label{example:1}
Let $\theta \in \Theta$ with $N = \{1,2\}$, $\lambda = (1, 5)$, $C = (1,1)$, $h=(0,0)$, $\mu =  (1, 1)$, and  $d = (4, 1)$. Under full pooling, it is optimal to stock one spare part for both DoDs individually, leading to  $c^{\theta}_F(\{1\})  = 2$ and $c^{\theta}_F(\{2\}) = 4 \frac{1}{6}$. If the DoDs collaborate, it is optimal to stock two spare parts, leading to $c^{\theta^F}(\{1,2\}) = 6 \frac{3}{13}$. Note that $c^{\theta}_F(\{1\}) + c^{\theta}_F(\{2\}) < c^{\theta}_F(\{1,2\})$, implying that the core is empty.
\end{example}

In Example \ref{example:1} we see that the second player, the one with lower downtime costs, is actually overdemanding from the pool, while the first player is not willing to pay for that
in the end. With threshold pooling, the second player will be limited to overdemand from the pool, as we demonstrate in the upcoming example.

\begin{example}
Reconsider the spare parts pooling situation of Example \ref{example:1}, but now DoDs apply threshold pooling. For the individual DoDs, it is still optimal to stock one spare part, leading to $c^{\theta}(\{1\})  = 2$ and $c^{\theta}(\{2\}) = 4 \frac{1}{6}$. If the DoDs collaborate, it is optimal to stock two spare parts, with DoD 1 having the freedom to  always demand, while DoD 2 can only demand if two parts are on stock. The associated costs are $c^{\theta}(\{1,2\} = 5 \frac{2}{11}$. Note that $c^{\theta}({1}) + c^{\theta}({2})  \geq   c^{\theta}({1,2})$, implying that the core is non-empty.
\end{example}

\section{Core non-emptiness of threshold pooling games}\label{sec:mainresult}

In the previous section, we demonstrated that the core of full-pooling games can be empty. This is, however, not possible if players apply threshold pooling. In this section, we prove this via a two step procedure. In the first step, we show that there exists a one-dimensional value function whose output coincide with the original multi-dimensional value function of Lemma \ref{grotemdp} for some relevant states. Second, we use the one-dimensional value function to prove that our threshold pooling game is balanced, implying core non-emptiness. 


\subsection{Step 1: Reduction to one-dimensional space}

We first present some structural results for value function $W^S_t$ of coalition $S \subseteq N$ for all $t \in \mathbb{N}_0$. These results tell us that when storing a repaired spare part, an optimal policy would select a player with lowest holding costs among the players with storage capacity available (Lemma \ref{voorraadgoedkooppakken1}). Conversely, when satisfying demand, an optimal policy would withdraw a spare part from a player with highest holding costs among the players with spare parts available (Corollary \ref{voorraadgoedkooppakken2}).

\begin{theorem5} Let $\theta \in \Theta$ and $S \subseteq N$. For all $t \in \mathbb{N}_0$ it holds that
$$W_t^S(z+e^S(q)) \leq W^S_t(z+e^S(w))$$

\noindent for all $q,w \in S$ with  $q \leq w$ and all $z \in \mathscr{Z}^S$ for which $z+e^S(q), z+e^S(w) \in \mathscr{Z}^S$. \label{voorraadgoedkooppakken1}
\end{theorem5}
A similar result follows by applying Lemma~\ref{voorraadgoedkooppakken1} to $z'= z-e^S(w)-e^S(q)$.
\begin{theorem4} Let $\theta \in \Theta$ and $S \subseteq N$. For all $t \in \mathbb{N}_0$ it holds that
$$W_t^S(z-e^S(q)) \geq W^S_t(z-e^S(w))$$

\noindent for all $q,w \in S$ with  $q \leq w$ and all $z \in \mathscr{Z}^S$ for which $z-e^S(q), z-e^S(w) \in \mathscr{Z}^S$. \label{voorraadgoedkooppakken2} 
\end{theorem4}

The results of Lemma \ref{voorraadgoedkooppakken1} and Corollary \ref{voorraadgoedkooppakken2} tell us that, in the long-run, it is optimal for the authority to stock spare parts at the cheapest holding cost locations and to demand from the most expensive holding cost locations. Put differently, in the long-run, we are only interested in states of $\mathscr{Z}^S$ with $S \subseteq N$ for which there exists a $j \in S$ such that for all locations $i \in N$ with $i < j$  we have $z_i = C_i$, for location $j$ we have $0 \leq z_{j} < C_i$  and for all locations $i \in N$ with $i > j$ we have $z_i = 0$. We formalize this state space below. 

\begin{theorem1} Let $\theta \in \Theta$. Then, for all $S \subseteq N$, we define
\begin{equation*}\overline{\mathscr{Z}}^{S} = \{z \in \mathscr{Z}^S \hspace{1mm} \vert \hspace{1mm} \exists \hspace{1mm} j \in S: z_i=C_i \hspace{1mm}\mbox{ for } i<j \mbox{ and } z_i =0 \mbox{ for } i>j\} \end{equation*} \end{theorem1}

It is important to realize that for this new state space $\overline{\mathscr{Z}}^S$, there is actually no reason to keep track of the number of spare parts per local stock point. That is, due to the specific structure of assigning spare parts to demand and storing repaired spare parts to locations, it suffices to keep track of the total number of spare parts in stock only. 
To formalize this, let us introduce a new state space $\mathscr{Y}^S = \{0,\ldots,C_S\}$ with $C_S = \sum_{i \in S} C_i$ for all $S \subseteq N$.

Next, we introduce a function $H^{S}$ for all $S \subseteq N$, which represents the total holding costs for a given total number spare parts on stock, a coalition $S \subseteq N$ and associated capacities $(C_i)_{i \in S}$.

\begin{theorem1} 
\label{def:holdingcosts}Let $\theta \in \Theta$, $S \subseteq N$
. For all $y \in \{0,\ldots, 
 C_S\}$, define
$$H^{S}(y) = \sum_{i \in S} h_i^* \cdot \max\{0, \min\{y - \sum_{j \in S: j<i}C_j,C_i\}\}.$$ \end{theorem1}
Finally, we present the one-dimensional value function.

\begin{theorem1} \label{valuekie} Let $\theta \in \Theta$ and $S \subseteq N$. For all $y \in \mathscr{Y}^S$ and all $t \in \mathbb{N}_0$, define
\begin{equation*} \begin{aligned}
 V^S_{t+1}(y) &=  \sum_{i \in S} \left[\lambda^*_i \min_{l \in \{0,\min\{y,1\}\}} \left\{V^S_t(y-l) + (1-l)d_i \right\} + \mu^*_i\min_{l \in \left\{0, \min{\left\{1,C_S-y\right\}}\right\}} V^S_t(y+l) \right] \\  &  \hspace{2mm}+ H^{S}(y) + \left(1- \sum_{i \in S} \left[ \lambda^*_i + \mu^*_i\right] \right) \cdot V^S_t(y)\end{aligned} \end{equation*}
 \noindent and $V^S_0(y) = 0$ for all $y \in \mathscr{Y}^S$. \smallskip
\end{theorem1}

We now prove that $W_t^S$ and $V_t^S$ coincide for all $t \in \mathbb{N}_0$ and all relevant states $z \in \overline{\mathscr{Z}}^S  \subseteq \mathscr{Z}$. 
\begin{theorem5} 
\label{lem:cardinality}
Let $\theta \in \Theta$ and $S \subseteq N$. Then, for all $z \in \overline{\mathscr{Z}}^S$ it holds that
  \begin{equation*} W_t^S(z) = V_t^S(\vert \vert z \vert \vert_{_1}) \mbox{ for all } t \in \mathbb{N}_0.\end{equation*}
\end{theorem5}


Because Lemma 2 holds for any $z \in \mathscr{Z}^S$, we can use Lemma \ref{lem:cardinality} to conclude that for any $S \subseteq N \backslash \{\emptyset\}$,
\begin{equation*}c^{\theta}(S) = \gamma \cdot \lim_{t \to \infty} \frac{V_t^S(\vert \vert z \vert \vert_{_1})}{t} \mbox{ for all } z \in \overline{\mathscr{Z}}^S,\end{equation*}
and consequently, because $(C_i)_{i \in S} \in \overline{\mathscr{Z}}^S$, we also have 
\begin{equation}\label{eq:useful}c^{\theta}(S) = \gamma \cdot \lim_{t \to \infty} \frac{V_t^S\left(C_S\right)}{t}.\end{equation}

\subsection{Step 2: Showing balancedness}
Take any $\theta \in \Theta$ and let $\mathscr{B} \subseteq 2^N \backslash \{\emptyset\}$ be a minimal balanced collection. Using the notation and results of Section~\ref{sec:balancedness}, we will show that $\sum_{S \in \mathscr{B}} b_S \cdot c^{\theta}(S) $ is greater than or equal to $\alpha \cdot c^{\theta}(N)$, which by Theorem~\ref{thm:degrotestelling} proves that the core of our game is non-empty. As a first step in our proof, we will construct $b_S$ labeled \emph{copies} for all coalitions $S$ in $\mathscr{B}$. The collection of these copies of the various coalitions is defined next.
\begin{theorem1} Let $\theta \in \Theta$ and $ \mathscr{B} \subseteq 2^N \backslash \{\emptyset\}$ be a minimal balanced collection. Then, we define
 \begin{equation*} \begin{aligned} \mathscr{L} &= \bigg\{ \hspace{2mm} (S,k) \hspace{2mm} \vert \hspace{2mm} S \in \mathscr{B}, k \in \{1,2,\ldots,b_S\}\bigg\}. \end{aligned} \end{equation*}\end{theorem1}

We will refer to these labeled copies as \emph{labeled coalitions}. For each labeled coalition $(S, k) \in \mathscr{L}$, we denote $V^{S,k}=V^S$, $C_{S,k}=C_S$, $H^{S,k} = H^S$, and $\mathscr{Y}^{S,k} = \mathscr{Y}^S$. 
By construction of $\mathscr{L}$ and equation (\ref{eq:useful})
, we can now express $\sum_{S \in \mathscr{B}} b_S \cdot c^{\theta}(S)$ in terms of the new value functions $V^{S,k}$.
 \begin{theorem5} \label{eerstestap} For every $\theta \in \Theta$ it holds for any minimal balanced collection $\mathscr{B} \subseteq 2^N \backslash \{\emptyset\}$ that
\begin{equation} \sum_{S \in \mathscr{B}} b_S \cdot c^{\theta}(S) = \gamma \cdot \lim_{t \to \infty} \frac{1}{t} \cdot \sum_{S \in \mathscr{B}} \sum_{k=1}^{b_S} V^{S,k}_t(C_{S,k}).\label{eq:balancedexpression}\end{equation}
 \end{theorem5}

The next step is to construct a combined value function (of some yet unspecified MDP) with a state space that keeps track of the inventory level of every labeled coalition $(S, k) \in \mathscr{L}$, an action space that consists of all possible actions per labeled coalition $(S, k) \in \mathscr{L}$ given its
inventory level, and for which the related costs coincide with $\sum_{S \in \mathscr{B}} \sum_{k=1}^{b_S} V^{S,k}_t(C_{S,k})$ for
all $t \in \mathbb{N} \cup \{0\}$. We start with introducing the new state space. For notational convenience, we will (often) use variable $z$ to represent a labeled coalition $(S,k)$.

\begin{theorem1} Let $\theta \in \Theta$ and $ \mathscr{B} \subseteq 2^N \backslash \{\emptyset\}$ be a minimal balanced collection. Then, we define
 \begin{equation*} \begin{aligned} \mathscr{Y}^{\mathscr{B}} &= \left\{ \hspace{3mm} (r_{z})_{z \in \mathscr{L}} \hspace{3mm} \bigg\vert \hspace{2mm} r_{z} \in \left\{0,1,\ldots, C_{z}\right\} \hspace{2mm} \forall z \in \mathscr{L} \right\}. \end{aligned} \end{equation*}\end{theorem1}

Subsequently, we introduce the new action spaces and value function.

\begin{theorem1} Let $\theta \in \Theta$ and $\mathscr{B} \subseteq 2^N \backslash \{\emptyset\}$ be a minimal balanced collection. Then, for all $r=(r_z)_{z \in \mathscr{L}} \in \mathscr{Y}^{\mathscr{B}}$ and all $i \in N$ we define
  \begin{equation*} \begin{aligned} \mathscr{A}^{\mathscr{B}}_{i,-}(r) &= \left\{ \hspace{2mm} (l_z)_{z \in \mathscr{L}} \hspace{2mm} \bigg\vert \hspace{2mm} \begin{matrix} & l_{S,k}=l_z \in \{0,\min\{1,r_z\}\} \hspace{2mm}&&\mbox{ if } i \in S \\
&l_{S,k}=l_z=0 &&\mbox{ if } i \not \in S \end{matrix} \right\}; \\
  \mathscr{A}^{\mathscr{B}}_{i,+}(r) &= \left\{ \hspace{2mm} (l_z)_{z \in \mathscr{L}} \hspace{2mm} \bigg\vert \hspace{2mm} \begin{matrix} &l_{S,k}=l_z \in \{0,\min\{1,C_z - r_z\}\} \hspace{2mm}&&\mbox{ if } i \in S \\
&l_{S,k}=l_z = 0 && \mbox{ if } i \not \in S \end{matrix} \right\}.
\end{aligned} \end{equation*}\label{def:actionpercoalition}
\end{theorem1}

\begin{theorem1} \label{tweedeval} Let $\theta \in \Theta$ and $\mathscr{B} \subseteq 2^{N} \backslash\{\emptyset\}$ be a minimal balanced collection. Then, for all $r \in \mathscr{Y}^{\mathscr{B}}$ and all $t \in \mathbb{N}_0$, we define the value function as
\begin{equation*} \begin{aligned} V^{\mathscr{B}}_{t+1}(r) =  &\sum_{i \in N} \left[\lambda^*_i \min_{l \in \mathscr{A}^{\mathscr{B}}_{i,-}(r)}\bigg\{ (\alpha - \vert \vert l \vert \vert_{_1}) d_i + V^{\mathscr{B}}_t(r-l) \bigg\} \right. \\
&\hspace{10mm}+ \left. \mu^*_i \min_{l \in \mathscr{A}^{\mathscr{B}}_{i,+}(r)} \left\{ V^{\mathscr{B}}_t(r+l)  \right\} \right] + \sum_{S \in \mathscr{B}} \sum_{k=1}^{b_S} H^{S,k}(r_{S,k}). 
  \end{aligned} \end{equation*}
  \noindent with $V^{\mathscr{B}}_{0}(r) = 0$ for all $r \in \mathscr{Y}^{\mathscr{B}}$.
  \end{theorem1}
 

The new value function can be interpreted in the following way. Let $i \in N$.
With probability $\lambda_i^*$ there is a demand for all labeled coalitions $(S, k) \in \mathscr{L}$ for which
$i \in S$. Each such labeled coalition $(S, k)$ has, except for the case with $r_{S,k} = 0$, the possibility to accept demand ($l_{S,k} = 1$), and always the possibility to reject demand ($l_{S,k} = 0$). For all (other) labeled coalitions $(S, k) \in \mathscr{L}$ for which $i \not \in S$ it holds that there is no demand arrival and so $l_{S,k} = 0$. Based on these decisions, total costs equal $(\alpha - \vert \vert  l \vert \vert_1 ) d_i$ and one transits to state $r - l$. With probability $\mu_i^*$ there is a spare part repair for each labeled  coalition $(S, k) \in \mathscr{L}$ for which $i \in S$. For each labeled coalition $(S, k) \in \mathscr{L}$ with $i \in S$ and $r_{S,k} < C_{S,k}$ one can accept ($l_{S,k} = 1$) or reject ($l_{S,k} = 0$) the spare part. An accept will follow if the current number of spare parts on stock is smaller than what this labeled coalition would stock optimally and a reject will follow otherwise. 
For all (other) labeled coalitions $(S, k) \in \mathscr{L}$ for which $i \not \in S$ and for which $i \in S$ with $r_{S,k} = C_{S,k}$ it holds that there is no spare part
repair and so $l_{S,k} = 0$. Based on the decisions made, one subsequently transits to state $r + l$.

  We now prove that the new value function and $\sum_{S \in \mathscr{B}} \sum_{k=1}^{b_S} V^{S,k}_t(C_{S,k})$ coincide for all $t \in \mathbb{N}_+$.

\begin{theorem5} \label{combinatie4} Let $\theta \in \Theta$ and $\mathscr{B} \subseteq 2^{N} \backslash \{\emptyset\}$ be a minimal balanced collection. Then, for all $r \in \mathscr{Y}^{\mathscr{B}}$ and all $t \in \mathbb{N}_0$ it holds that
\begin{equation*} \sum_{S \in \mathscr{B}}\sum_{k=1}^{b_S} V^{S,k}_t(r_{S,k}) = V^{\mathscr{B}}_t(r).\end{equation*} \end{theorem5}

Note that the action space of 
$V_t^{\mathscr{B}}$ 
is restricted. For instance, upon a
demand arrival of player $i \in S$, it is not possible to accept a single demand for labeled
coalition $(S, k) \in \mathscr{L}$ for which $i \in S$ and $r_{S,k} = 0$, while other labeled coalitions
$(S, k) \in \mathscr{L}$ may still be able to accept it. A similar reasoning holds
for repair completion. It is not possible to store a spare part for any labeled coalition
$(S, k) \in \mathscr{L}$ for which $r_{S,k} = C_{S,k}$, while other not fully replenished labeled
coalitions may be able to store the spare part. The next step is to introduce a value function that incorporates these
extended possibilities. So, we introduce a value function (related to some unspecified
MDP), that coincides with the value function of Definition \ref{tweedeval}, except for a relaxed action
space. In order to do so, we first need to introduce this relaxed action space.

\begin{theorem1} Let $\theta \in \Theta$ and $\mathscr{B} \subseteq 2^{N} \backslash \{\emptyset\}$ be a minimal balanced collection. Then, for all $r \in \mathscr{Y}^{\mathscr{B}}$ and all $i \in N$ we define
  \begin{equation*} \begin{aligned}   \widehat{\mathscr{A}}^{\mathscr{\hspace{1.3mm} B}}_{i,-}(r) &= \left\{ \hspace{2mm} (l_z)_{z \in \mathscr{L}} \hspace{2mm} \bigg\vert \hspace{2mm}\forall z \in \mathscr{L}:  l_z \in \{0,1,\ldots,\alpha\} \hspace{2mm} , \hspace{2mm} \sum_{z \in \mathscr{L}} l_z\leq \alpha, \hspace{2mm} r - l \in \mathscr{Y}^{\mathscr{B}} \right\}; \\
\widehat{\mathscr{A}}^{\hspace{1.3mm}\mathscr{B}}_{i,+}(r) &= \left\{ \hspace{2mm} (l_z)_{z \in \mathscr{L}} \hspace{2mm} \bigg\vert \hspace{2mm}  \forall z \in \mathscr{L}:  l_z \in \{0,1,\ldots,\alpha\}, \hspace{2mm} \sum_{z \in \mathscr{L}} l_z\leq \alpha, \hspace{2mm} r + l \in \mathscr{Y}^{\mathscr{B}} \right\}. \\
\end{aligned} \end{equation*}\label{def:actionspace2}
\end{theorem1} 
Please, observe that the action spaces of Definition \ref{def:actionspace2} contain the action spaces of Definition \ref{def:actionpercoalition}, as formalized in the following lemma (the proof is
straightforward and for this reason omitted).

\begin{theorem5} \label{kakieplant} Let $\theta \in \Theta$ and $\mathscr{B} \subseteq 2^{N} \backslash \{ \emptyset\}$ be a minimal balanced collection. Then, for all $r \in \mathscr{Y}^{\mathscr{B}}$ and all $i \in N$ it holds that $\mathscr{A}^{\mathscr{B}}_{i,-}(r) \subseteq \widehat{\mathscr{A}}^{\hspace{1.3mm}\mathscr{B}}_{i,-}(r)$ and $\mathscr{A}^{\mathscr{B}}_{i,+}(r) \subseteq \widehat{\mathscr{A}}^{\hspace{1mm} \mathscr{B}}_{i,+}(r)$.
\end{theorem5}

By relaxing the action space, we also allow spare parts to be stored differently, i.e., stored in a cheaper way. Instead of storing spare parts per labeled coalition only, it is now possible to stock them anywhere, taking into account that each local stock point $i \in N$ has a storage capacity of $\alpha \cdot C_i$. The holding costs associated to an optimal storage of spare parts is denoted by 
$$H^{N}_{\alpha}(y) = \sum_{i \in N} h_i^* \cdot \max\{0, \min\{y - \sum_{j \in N: j<i}\alpha \cdot C_j,\alpha \cdot C_i\}\} \mbox{ for all } y \in \{0,\ldots, 
 \alpha \cdot C_N\}.$$
We are now ready to define the value function using the  new action spaces and holding costs. 
\begin{theorem1} \label{valuefunctionsint} Let $\theta \in \Theta$ and $\mathscr{B} \subseteq 2^{N} \backslash \{\emptyset\}$ be a minimal balanced collection. For all $r \in \mathscr{Y}^{\mathscr{B}}$ and all $t \in \mathbb{N}_0$, we define
\begin{equation*} \begin{aligned}\widehat{V}_{t+1}^{\mathscr{B}}(r) =  &\sum_{i \in N} \left[\lambda^*_i \min_{l \in \widehat{\mathscr{A}}^{\hspace{1.3mm} \mathscr{B}}_{i,-}(r)}\bigg\{ (\alpha - \vert \vert l \vert \vert_{_1}) d_i + \widehat{V}_t^{\mathscr{B}}(r-l) \bigg\} \right. \\
 & \hspace{10mm} + \left. \mu^*_i \min_{l \in \widehat{\mathscr{A}}^{\hspace{1.3mm} \mathscr{B}}_{i,+}(r)} \left\{ \widehat{V}_t^{\mathscr{B}}(r+l)  \right\} \right] + H^{N}_{\alpha}(\vert \vert r \vert \vert_{_1}).
  \end{aligned} \end{equation*}

  \noindent with $\widehat{V}^{\mathscr{B}}_{0}(r) = 0$ for all $r \in \mathscr{Y}^{\mathscr{B}}$. \end{theorem1}
 
 We now formalize that the value function 
 of Definition~\ref{valuefunctionsint} relaxes the value function 
 of Definition~\ref{tweedeval}.
\begin{theorem5} \label{tweedestap} Let $\theta \in \Theta$ and $\mathscr{B} \subseteq 2^{N} \backslash \{\emptyset\}$ be a minimal balanced collection. For all $r \in \mathscr{Y}^{\mathscr{B}}$ and all $ t \in \mathbb{N}_0$ it holds that
\begin{equation*} V^{\mathscr{B}}_{t}(r) \geq \widehat{V}_{t}^{\mathscr{B}}(r). \end{equation*} \end{theorem5}

We will next show that 
$\widehat{V}_t^{\mathscr{B}}$
is the same for all states $r,r' \in \mathscr{Y}^{\mathscr{B}}$ 
 for
which $\vert \vert r\vert \vert_1 = \vert \vert r' \vert \vert_1$ and so the decisions made in these states exhibit a
similar equivalence. Consequently, it suffices to focus on an alternative MDP whose states and actions depend on the total number of spare parts only. We introduce this MDP next, before formalizing the correspondence in Lemma \ref{theorem:bnaara}.

\begin{theorem1} \label{valuesint2} Let $\theta \in \Theta$ and $\mathscr{B} \subseteq 2^{N} \backslash \{\emptyset\}$ be a minimal balanced collection. Then, for all $j \in \{0,1,\ldots, \alpha \cdot C_N\}$ and all $t \in \mathbb{N}_0$ we define
\begin{equation*} \begin{aligned}V^{\alpha}_{t+1}(j) &= \sum_{i \in N} \left[\lambda^*_i\min_{l \in \{0,.., \min\{\alpha,j\}\}}\bigg\{ (\alpha - l) d_i + V^{\alpha}_t(j-l) \bigg\} \right. \\
&\hspace{15mm} \left. +  \mu^*_i \min_{l \in \{0,.., \min\{\alpha, \alpha \cdot C_N - j\}\}}  V^{\alpha}_t(j+l) \right] + H^{N, \alpha}(j)
  \end{aligned} \end{equation*}
  \noindent with $V^{\alpha}_0(j) = 0$ for all $j \in \{0,1,\ldots,\alpha \cdot C_N\}$.
  \end{theorem1}

\begin{theorem5} \label{theorem:bnaara} Let $\theta \in \Theta$ and $\mathscr{B} \subseteq 2^{N} \backslash \{\emptyset\}$ be a minimal balanced collection. For all $r \in \mathscr{Y}^{\mathscr{B}}$ we have \label{ano}
$$\widehat{V}^{\mathscr{B}}_t(r) = V^{\alpha}_t( \vert\vert r\vert\vert_{_1}) \hspace{7mm} \text{ for all } t \in \mathbb{N}_0.$$
\end{theorem5}

In $V^{\alpha}_t$
demand arrives in batches of $\alpha$ and\, similarly, batches of $\alpha$ parts arrive whenever a repair completes. Moreover, the total inventory capacity reads $\alpha C_N$. The next objective is to derive a correspondence between $V^{\alpha}_t$, and a related MDP where batch sizes are $1$ and the total inventory capacity is scaled down by factor $\alpha$, i.e., the capacity is $C_N$. To this end, we investigate two properties. 

\begin{theorem5} \label{allesofniks}
  Let $\theta \in \Theta$ and $\mathscr{B} \subseteq 2^N \backslash \{\emptyset\}$ be a minimal balanced collection. For all $t \in \mathbb{N}_0$ it holds that \medskip

 \hspace{4mm} \noindent $(i)$ $V^{\alpha}_t(j) + V^{\alpha}_t(j+2) \geq 2\cdot V^{\alpha}_t(j+1)$ for all $j \in \{0,1,\ldots,\alpha \cdot C_N - 2\}$; \smallskip

 \hspace{3mm} \noindent $(ii)$ $V^{\alpha}_t(k + j) + V^{\alpha}_t(k  + j+2) = 2\cdot V^{\alpha}_t(k + j +1)$ for all $j \in \{0,1,\ldots,\alpha -2\}$ \smallskip

 \hspace{12mm} and all $k \in \{0,\alpha, 2\alpha,\ldots,(C_N-1)\alpha\}$.
\end{theorem5}
The lemma demonstrates that the value function is convex, and piecewise linear. As a consequence, when a batch of $\alpha$ demands arrives, it is either satisfied or rejected in full. Moreover, a repair completion either corresponds to an increase in the inventory of $\alpha$ parts, or no increase at all. 
From this, we can conclude that the states of $V^{\alpha}_t$
that are multiples of $\alpha$ depend on the states of $V^{\alpha}_t$ that are multiples of $\alpha$ only. This allows us to rewrite value function $V^{\alpha}_t$
into $\alpha$ - times value function $V^N_t$.

  \begin{theorem5} \label{theorem:alfadelen} Let $\theta \in \Theta$ and $\mathscr{B} \subseteq 2^N \backslash \{\emptyset\}$ be a minimal balanced collection. For all $j \in \{0,\alpha,\ldots,C_N \cdot \alpha\}$ and all $t \in \mathbb{N}_0$ it holds that\label{alphading}
  \begin{equation*} V^{\alpha}_t(j) = \alpha \cdot V^N_t\left(\frac{j}{\alpha}\right)\end{equation*}
  \end{theorem5}
With this link established, we have everything in place to prove our main result.
\begin{theorem0} Threshold pooling games are balanced. \label{hoofdstelling} \end{theorem0}

\noindent \emph{\textbf{Proof :}} Let $\theta \in \Theta$ and $(N,c^{\theta})$ be the associated threshold pooling game and $\mathscr{B} \subseteq 2^{N} \backslash \{\emptyset\}$ be a minimal balanced collection. Then, observe that
\begin{equation*} \begin{aligned} \sum_{S \in \mathscr{B}} b_S \cdot c^{\theta}(S) &= \gamma \cdot \lim_{t \to \infty} \frac{1}{t} \sum_{S \in \mathscr{B}} \sum_{k=1}^{b_S} V^{S,k}_t(C_{S,k}) \geq \gamma \cdot \lim_{t \to \infty} \frac{ \widehat{V}^{\mathscr{B}}_{t}((C_{S,k})_{(S,k) \in \mathscr{L}}) }{t} \\
&=\gamma \cdot  \lim_{t \to \infty} \frac{ V^{\alpha}_{t}(\alpha \cdot C_N)}{t}  =  \gamma \cdot \lim_{t \to \infty} \alpha \cdot \frac{V^N_t(C_N)}{t} = \alpha \cdot c^{\theta}(N) \\
\end{aligned} \end{equation*} 
\noindent The first equality holds by Lemma \ref{eerstestap}. The inequality holds by Lemma \ref{combinatie4} and Lemma
\ref{tweedestap}. The second equality holds by Lemma \ref{theorem:bnaara} and the fact that $\sum_{(S,k) \in \mathscr{L}}C_{S,k} = \alpha \cdot C_N$. The
third equality holds by Lemma \ref{theorem:alfadelen}. 
The last equality holds by taking $\alpha$ outside
the limit (which is allowed as it is a constant) and subsequently using the relationship in equation (\ref{eq:useful}). 
Finally, note that this establishes \eqref{balanced} for any threshold pooling game, which implies that such games are balanced in general. This concludes the proof. \bigskip $\hfill \square$

 Based on Theorem \ref{thm:degrotestelling}, the next result follows immediately.

\begin{theorem4} \label{thm:degrotestelling}Threshold pooling games have a non-empty core.\end{theorem4}

\section{Concluding remarks}\label{sec:conclusion}

 Our results have demonstrated that the type of pooling strategy plays a crucial role in the success of a joint spare parts pool. More precisely, we showed that a joint spare parts pool may not last long --or even not arise-- if full pooling is applied, while a joint spare parts pool is possible under threshold pooling. We did so by studying two cooperative games, one for full pooling and one for threshold pooling. For both games, the coalitional values were formulated as the  long-term average costs of an underlying MDP. 
 
The results in this study are established under some assumptions. One of them is that 
 DoDs have a maximal capacity per local stock point. In practice, however, these capacities might not play a dominant role. We would like to stress that our games also apply to such a setting, namely by iteratively increasing the capacities until they are not binding anymore for any coalition. 
Next, we also decided to not explicitly model travel times (and associated downtime costs) when spare parts are transported between local stock points. This assumption is reasonable when stock points are located within 1 or 2 hours travel time for each other, but not when, for instance, stock points are located at different continents. For such a setting, it could be that core non-emptiness is no longer guaranteed --even under threshold pooling.  
An interesting future research direction for such a setting would then be to study the weak epsilon core (see, e.g., \cite{shapley1966quasi}). This epsilon could, for instance, represent a subsidy paid by a third-party (e.g., NATO or the European Defence Agency) to each DoD such that everyone is still willing to join the spare parts pool.
Finally, we considered in this paper a setting with \emph{one} specific repairable part only. In practice, however, DoDs might be interested in pooling  multiple spare parts simultaneously. An interesting future direction would therefore be to study how core non-emptiness is preserved if multiple spare parts are pooled simultaneously. We would like to stress that our results cannot be applied to this setting directly. Studying a threshold pooling for each type of spare part individually could, for instance, lead to an overestimation of downtime costs.

\bibliographystyle{abbrvnat}
\bibliography{reference}

\newpage 

\section{Appendix}

\noindent \underline{Proof of Lemma 1} \bigskip

\noindent Let $z \in  \mathscr{Z}^S$ and $t \in \mathbb{N}_0$. Then
\begin{equation*} \begin{aligned} W^S_{t+1}(z) =& \min_{a \in \mathscr{A}^S(z)} \left\{ C^S(z,a) + \sum_{z' \in \mathscr{Z}^S} p(z' \vert z,a) \cdot W^S_t(z')\right\} \\
=& \min_{a \in \mathscr{A}^S(z)} \left\{ \sum_{i \in S} \lambda^*_i  \left(1-\sum_{j \in S:z_j>0}a_{ij}^{-}\right)d_i + \sum_{i \in S
} \lambda^*_i \sum_{j \in S: z_j > 0} a_{ij}^{-}  W^S_t(z - e^S(j)) \right. \\
& \left.  +\sum_{i \in S
} \mu^*_i \sum_{j \in S:z_j<C_j} a_{ij}^{+}  W^S_t(z+e^S(j))\right\}  \\
&+ \left( 1 - \sum_{i \in S}\left[ \lambda^*_i \sum_{j \in S: z_j >0}a_{ij}^{-} +\mu^*_i \sum_{j \in S: z_j < C_j} a_{ij}^{+} \right]\right) W^S_t(z) + \sum_{i \in S} h^*_i \cdot z_i \\
=& \min_{a \in \mathscr{A}^S(z)} \left\{\sum_{i \in S} \left[\lambda^*_i  \left( \left(1-\sum_{j \in S:z_j>0}a_{ij}^{-}\right)\left(d_i+ W^S_t(z)\right) + \sum_{j \in S:z_j>0} a_{ij}^{-} W^S_t(z-e^S(j)) \right)  \right.\right. \\
&  \left. + \mu^*_i \sum_{j \in S:z_j<C_j} a_{ij}^{+} W^S_t(z+e^S(j)) \bigg]
+ \mu^*_i \left( 1 - \sum_{j \in S:z_j<C_j} a_{ij}^{+} \right) W^S_t(z) \bigg]  \right. \\
&  \left.  + \left(1 - \sum_{i \in S} \left[ \lambda^*_i + \mu^*_i\right] \right) W^S_t(z)\right\} + \sum_{i \in S} h_i^* \cdot z_i\\
=&\sum_{i \in S} \left[\lambda^*_i \min_{x \in \mathscr{A}^S_{-}(z)} \left\{(1-\vert \vert x \vert \vert_{_1})d_i + W^S_t(z-x) \right\} + \mu_i^*\min_{x \in \mathscr{A}^S_{+}(z)} W^S_t(z+x) \right]  \\
 &  \hspace{17mm}+ \left(1- \sum_{i \in S} \left[ \lambda^*_i + \mu^*_i\right] \right) \cdot W^S_t(z) + \sum_{i \in S} h_i^* \cdot z_i.
\end{aligned} \end{equation*}

\noindent The first equality holds by (\ref{valuefunctiondef}). In the second equality, we use the definition of $C^S(z,a)$ and $p(z' \vert z,a)$, in combination with the fact that $\sum_{j \in S}a_{ij}^{+} = \sum_{j \in S: z_j>0} a_{ij}^{+}$ and $\sum_{j \in S}a_{ij}^{-} = \sum_{j \in S: z_j>0} a_{ij}^{-}$ for all $i \in S$ and all $z \in \mathscr{Z}^S$. Moreover, we use that $\sum_{i \in S}h_i^* \cdot z_i$ is independent of the action set. In the third equality, we did some rewriting. The last equality holds since a minimum of a sum of independent terms equals the sum of all these individual terms evaluated at their minimum. This concludes the proof. $\hfill \square$ \bigskip

\noindent \underline{Proof of Lemma 2} \bigskip

\noindent The first equality holds by uniformization, which is allowed if transition rates are bounded and the MDP is multichain (see \citet[p.568]{puterman2014markov}). Notice that interarrival times of demands as well as repair times are exponentially distributed with rates that are bounded from above. In addition, for every stationary policy, there exist one or multiple recurrent classes. Hence, the MDP is multichain. With respect to the second equality, observe that state space $\mathscr{Z}^S$ and action space $\mathscr{A}^S(z)$ for all $z \in \mathscr{Z}^S$ of the MDP are finite. In addition, under stationary policy $f = (f_i(z))_{z \in \mathscr{Z}^S,i \in S}$ with $f_i(z) = (e^S(i),e^S(i))$ for all $i \in S$ and all $0<z<C_S$ and $f_i(z)=(0^S,e^S(i))$ for all $i \in S$ and all $z \in \mathscr{Z}^S$ for which $z_i=0$ and $f_i(z)=(e^S(i),0^S)$ for all $i \in S$ and all $z \in \mathscr{Z}^S$ for which $z_i=C_i$, every state 
$z \in \mathscr{Z}^S$ is accessible from any state 
$z' \in \mathscr{Z}^S$ after (possibly) some arrivals and some (one-by-one) production completions. Hence, the related Markov chain is irreducible. An irreducible Markov chain with finite state space is positive recurrent (see e.g., \citet{modica2012first}). Finally, observe that the long-run average costs per time epoch under policy $f$ are bounded (naturally) by $\sum_{i \in S} \lambda_i^* \cdot d_i + \sum_{i \in S} C_i \cdot h_i^*$ and as a result of \citet[Proposition 4.3]{sennott1996convergence}, the second equality follows. This concludes the proof. $\hfill \square$  \bigskip

\noindent \underline{Proof of Lemma 3} \bigskip

\noindent \textbf{Proof:} This proof is by induction. By definition of the value functions $W_0^S(z+e^S(q))=W_0^S(z+e^S(w))=0$ for all $q,w \in S$ with  $q \leq w$ and all $z \in \mathscr{Z}^S$ for which $z+e^S(q), z+e^S(w) \in \mathscr{Z}^S$. \bigskip

\noindent Let $t \in \mathbb{N}_0$ and assume that  $W_t^S(z+e^S(q)) \leq W^S_t(z+e^S(w))$ for all $q,w \in S$ with  $q \leq w$ and all $z \in \mathscr{Z}^S$ for which $z+e^S(q), z+e^S(w) \in \mathscr{Z}^S$. \bigskip

\noindent First, we show that for each $i \in N$, it holds that:
\begin{equation} \label{hulpstuk1} \begin{aligned} &\min_{x' \in \mathscr{A}^S_{-}(z+e^S(q))} \left\{W_t(z+e^S(q)-x') + (1-\vert \vert x' \vert \vert_{_{1}})d_i\right\} \\ \leq &\min_{x \in \mathscr{A}^S_{-}(z+e^S(w))} \left\{W_t(z+e^S(w)-x) +(1-\vert\vert x \vert\vert_{_{1}})d_i\right\} \end{aligned} \end{equation} 

\noindent for all $q,w \in S$ with  $q \leq w$ and all $z \in \mathscr{Z}^S$ for which $z+e^S(q), z+e^S(w) \in \mathscr{Z}^S$. \bigskip

\noindent Let $i \in N$. In addition, let $q,w \in S$ with  $q \leq w$ and $z \in \mathscr{Z}^S$ such that $z+e^S(q), z+e^S(w) \in \mathscr{Z}^S$. Let $x \in  \mathscr{A}^S_{-}(z+e^S(w))$. We distinguish between three cases. 
\bigskip

\noindent \emph{Case 1. $x =0^S$.} \medskip

\noindent Observe that $W_t(z+e^S(w)-x) + (1-\vert \vert x \vert \vert_{_{1}})d_i= W_t(z+e^S(w)) +d_i\geq W_t(z+e^S(q))$ $+ d_i = W_t(z+e^S(q)-x') +(1-\vert \vert x' \vert \vert_{_1})d_i$, with $x' = 0^S \in \mathscr{A}^S_{-}(z+e^S(q))$. The inequality holds by the induction hypothesis.    \bigskip

\noindent \emph{Case 2. $x=e^S(k)$ with $k \not = w$.} \bigskip

\noindent Observe that $W_t(z+e^S(w)-x) +(1-\vert\vert x \vert \vert_{_{1}})d_i = W_t(z+e^S(w)-e^S(k)) \geq W_t(z+e^S(q)-e^S(k)) = W_t(z+e^S(q)-x') + (1- \vert \vert x' \vert \vert_{_1})d_i$, with $x' =e^S(k) \in \mathscr{A}^S_{-}(z+e^S(q))$. The inequality holds by the induction hypothesis. \bigskip

\noindent \emph{Case 3. $x=e^S(w)$.} \bigskip

\noindent Observe that $W_t(z+e^S(w)-x) + (1 -\vert \vert x \vert \vert_{_{1}})d_i= W_t(z) = W_t(z +e^S(q) - x') + (1- \vert \vert x' \vert \vert_{_1})d_i$, with $x' =e^S(q) \in \mathscr{A}^S_{-}(z+e^S(q))$. \bigskip

\noindent Hence, we can conclude that (\ref{hulpstuk1}) is true. \bigskip

\noindent Secondly, we show  that for each $i \in N$, it holds that 
\begin{equation} \label{hulpstuk2} \min_{x' \in \mathscr{A}^S_{+}(z+e^S(q))} \left\{W_t(z+e^S(q)+x) \right\} \leq \min_{x \in \mathscr{A}^S_{+}(z+e^S(w))} \left\{W_t(z+e^S(w)+x) \right\} \end{equation}

\noindent for all $q,w \in S$ with  $q \leq w$ and all $z \in \mathscr{Z}^S$ for which $z+e^S(q), z+e^S(w) \in \mathscr{Z}^S$. \bigskip

\noindent Let $i \in N$. In addition, let $q,w \in S$ with  $q \leq w$ and $z \in \mathscr{Z}^S$ such that $z+e^S(q), z+e^S(w) \in \mathscr{Z}^S$. Let $x \in \mathscr{A}^S_{+}(z+e^S(w))$. We distinguish between three cases. \bigskip

\noindent \emph{Case 1.} $ x=0^S$. \bigskip

\noindent Observe that $W_t(z+e^S(w)+x) = W_t(z+e^S(w)) \geq W_t(z+e^S(q)) = W_t(z+e^S(q)+x')$, with $x' = 0^S \in \mathscr{A}^S_{+}(z+e^S(q))$. The inequality holds by the induction hypothesis.   \bigskip

\noindent \emph{Case 2.} $x =e^S(k)$ with $k \not=q$. \bigskip

\noindent Observe that $W_t(z+e^S(w)+x) = W_t(z+e^S(w)+e^S(k))\geq W_t(z+e^S(q)+e^S(k)) = W_t(z+e^S(q)+x')$, with $x' = e^S(k) \in \mathscr{A}^S_{+}(z+e^S(q))$. The inequality holds by the induction hypothesis.\bigskip

\noindent \emph{Case 3.} $x=e^S(q)$. \bigskip

\noindent Observe that $W_t(z+e^S(w)+x) = W_t(z+e^S(w)+e^S(q)) = W_t(z+e^S(q)+x')$, with $x' = e^S(w) \in \mathscr{A}^S_{+}(z+e^S(q))$. \bigskip

\noindent Hence, we can conclude that (\ref{hulpstuk2}) is true. \bigskip

\noindent Now, observe that
\begin{equation*} \begin{aligned} W_{t+1}(z + e^S(q)) =& \sum_{i \in S} \lambda_i^* \left[ \min_{x \in \mathscr{A}^S_{-}(z+e^S(q))} \left\{W_t(z+e^S(q)-x) + (1-\vert \vert x \vert \vert_{_{1}})d_i\right\} \right. \\
&\left. + \min_{x \in \mathscr{A}^S_{+}(z+e^S(q))} \left\{W_t(z+e^S(q)+x) \right\}  \right] + \sum_{i \in S} h_i^* \cdot z_i +  h_q \\
& + \left(1- \sum_{i \in S} \left[ \lambda^*_i + \mu^*_i\right] \right) \cdot W^S_t(z+e^S(q)) \\
\leq& \sum_{i \in S} \lambda_i^* \left[ \min_{x \in \mathscr{A}^S_{-}(z+e^S(w))} \left\{W_t(z+e^S(w)-x) + (1-\vert \vert x \vert \vert_{_{1}})d_i\right\} \right. \\
&\left. + \min_{x \in \mathscr{A}^S_{+}(z+e^S(w))} \left\{W_t(z+e^S(w)+x) \right\}  \right] + \sum_{i \in S} h_i^* \cdot z_i + h_q \\
&+ \left(1- \sum_{i \in S} \left[ \lambda^*_i + \mu^*_i\right] \right) \cdot W^S_t(z+e^S(q))
\end{aligned} \end{equation*}
\begin{equation*} \begin{aligned}
\hspace{30mm} \leq& \sum_{i \in S} \lambda_i^* \left[ \min_{x \in \mathscr{A}^S_{-}(z+e^S(w))} \left\{W_t(z+e^S(w)-x) + (1-\vert \vert x \vert \vert_{_{1}})d_i\right\} \right. \\
&\left. + \min_{x \in \mathscr{A}^S_{+}(y+e^S(w))} \left\{W_t(z+e^S(w)+x) \right\}  \right] + \sum_{i \in S} h_i^* \cdot z_i + h_w \\
&+ \left(1- \sum_{i \in S} \left[ \lambda^*_i + \mu^*_i\right] \right) \cdot W^S_t(z+e^S(w)) \\
& = W_{t+1}(z+e^S(w)),
\end{aligned} \end{equation*}
  
  \noindent where the first inequality holds by (\ref{hulpstuk1}) and (\ref{hulpstuk2}). The second inequality holds by the induction hypothesis and the relationship $h_q \leq h_w$, which holds by assumption as $q < w$. By the principle of mathematical induction, this completes the proof. $\hfill \square$  \bigskip

\noindent \underline{Proof of Lemma 4} \bigskip

\noindent \textbf{\emph{Proof:}}  This proof is by induction. By definition of the value functions $W_0^S(y)=V_0^S(\vert \vert y \vert \vert_{_1})=0$ for all $y \in \overline{\mathscr{Z}}^S$.\bigskip

\noindent Let $t \in \mathbb{N}_0$ and assume that  $W_t^S(y) = V^S_t(y)$ for all $y \in \overline{\mathscr{Z}}^S$. \bigskip

\noindent We define for all $z \in \overline{\mathscr{Z}}$
\begin{equation*} \begin{aligned} 
q^{+}(z) &= \left\{\begin{matrix} 0 & \mbox{ if } z=(C_i)_{i \in S}; \\ 
\min\{i \in S \vert z_i <C_i\} & \mbox{ else,} \end{matrix} \right.   \\
q^{-}(z) &= \left\{\begin{matrix} 0 & \mbox{ if } z=0^S; \\ 
\max\{i \in S \vert z_i >0\} & \mbox{ else.} \end{matrix} \right.
\end{aligned} \end{equation*}

\noindent Moreover, we extend the definition of $e^S$ of Section \ref{sec:MDP} by setting $e^S(0) = 0^S$. Note that $e^S(q^{+}(z)) \in \mathscr{A}^S_{+}(z)$  and $e^S(q^{-}(z)) \in \mathscr{A}^S_{-}(z)$ for all $z \in \overline{\mathscr{Z}}^S
$. \bigskip

\noindent Now, let $ z \in \overline{\mathscr{Z}}^S$. Observe that
\begin{equation*} \begin{aligned} W_{t+1}^S(z) =& \sum_{i \in S} \left[\lambda_i^* \min_{x \in \mathscr{A}^S_{-}(z)} \left\{ W_t(z-x) + (1 - \vert\vert x \vert\vert_{_{1}})d_i\right\} + \mu_i^* \min_{x \in \mathscr{A}^S_{+}(z)} \left\{W_t(z+x) \right\} \right]  \\
&+ \left(1- \sum_{i \in S} \left[ \lambda^*_i + \mu^*_i\right] \right) \cdot W^S_t(z) + \sum_{i \in S} h_i^* z_i \\
=&\sum_{i \in S} \left[\lambda_i^* \min_{x \in \{0^S,e^S(q^{-}(z))\}} \left\{ W_t(z-x) + (1 - \vert\vert x \vert\vert_{_{1}})d_i\right\} + \mu_i^* \min_{x \in \{0^S,e^S(q^{+}(z))\}} \left\{W_t(z+x) \right\} \right]  \\
&+ \left(1- \sum_{i \in S} \left[ \lambda^*_i + \mu^*_i\right] \right) \cdot W^S_t(z) + \sum_{i \in S} h_i^* \cdot z_i \\
=&\sum_{i \in S} \left[\lambda_i^* \min_{x \in \{0^S,e^S(q^{-}(z))\}} \left\{ V_t(\vert \vert z \vert \vert_{_1} - \vert \vert x \vert \vert_{_1}) + (1 - \vert \vert x \vert \vert_{_1})d_i\right\}  \right. \\
& \left. + \mu_i^* \min_{x \in \{0^S,e^S(q^{+}(z))\}} \left\{V_t(\vert \vert z \vert \vert_{_1} +\vert \vert x \vert \vert_{_1}) \right\} \right] + \left(1- \sum_{i \in S} \left[ \lambda^*_i + \mu^*_i\right] \right) \cdot V^S_t(\vert \vert z \vert \vert_{_1}) \\
& + \sum_{i \in S} h_i^* \cdot z_i  
 \end{aligned} \end{equation*}
\begin{equation*} \begin{aligned}
=&\sum_{i \in S} \left[\lambda_i^* \min_{k \in \{0,\min\{1,\vert \vert z \vert \vert_{_1}\}\}} \left\{ V_t(\vert \vert z \vert \vert_{_1} -k) + (1 - k)d_i\right\}  \right. \\
& \left. + \mu_i^* \min_{k \in \{0,\min\{1,C_s-\vert \vert z \vert \vert_{_1}\}\}} \left\{V_t(\vert \vert z \vert \vert_{_1} +k) \right\} \right] + \left(1- \sum_{i \in S} \left[ \lambda^*_i + \mu^*_i\right] \right) \cdot V^S_t(\vert \vert z \vert \vert_{_1}) \\
& + \sum_{i \in S} h_i^* \max\{0,\min\{C_i, \vert \vert z \vert \vert_{_1} - \sum_{j \in S: j < i}C_i\}\} \\
=& V_t^S(\vert \vert z \vert \vert_{_1}).
 \end{aligned} \end{equation*}

\noindent The first equality holds by Lemma \ref{grotemdp}. In the second equality, we use Lemma \ref{voorraadgoedkooppakken1} to conclude that an optimal decision $x \in \mathscr{A}^S_{-}(z)$ is either $0^S$ or $e^S(q^{-}(z))$ and use Corollary \ref{voorraadgoedkooppakken2} to conclude that an optimal decision $x \in \mathscr{A}^S_{+}(z)$ is either $0^S$ or $e^S(q^{+}(z))$. In the third equality, we use that for each $x \in \{0^S, e^S(q^{-}(z))\}$, we have $z-x \in \overline{\mathscr{Z}}^S$ and for each $x \in \{0^S, e^S(q^{+}(z))\}$, we have $z+x \in \overline{\mathscr{Z}}^S$, which allow us to apply the induction hypothesis three times. In the fourth equality, we introduce a new decision variable $k = \vert \vert z \vert \vert_{_1}$. In addition, we present the holding costs in a different way (which is possible, based on the specific structure of $z \in \overline{\mathscr{Z}}^S$). In the last equality, we use $H^S(\vert \vert z  \vert \vert_1) = \sum_{i \in S} h_i^* \max\{0,\min\{C_i, \vert \vert z \vert \vert_{_1} - \sum_{j \in S: j < i}C_i\}\}$ and subsequently use the definition of $V^S_t$. By the principle of mathematical induction, this completes the proof. \textcolor{white}{endproof} $\hfill \square$ \bigskip

\noindent \underline{Proof of Lemma 5} \bigskip

\noindent Let $\theta \in \Theta$ and $\mathscr{B} \subseteq 2^N \backslash \{\emptyset\}$ be a minimal balanced collection. It holds that
 $$ \begin{aligned} \sum_{S \in \mathscr{B}} b_S  \cdot c^{\theta}(S)  &= 
 \gamma  \cdot \sum_{S \in \mathscr{B}} \sum_{k=1}^{b_S} \lim_{t \to \infty} \frac{V^{S,k}_t(C_{S,k})}{t} \\
 &= \gamma \cdot \lim_{t \to \infty} \frac{1}{t}  \cdot  \sum_{S \in \mathscr{B}} \sum_{k=1}^{b_S} V^{S,k}_t(C_{S,k}). \end{aligned} $$
\noindent The first equality holds by exploiting all labeled coalitions and equation (\ref{eq:useful}). 
The second equality holds as all limits are well-defined and all sums are finite. 
$\hfill \square$ \bigskip

\noindent \underline{Proof of Lemma \ref{combinatie4}} \bigskip

\noindent \textbf{\emph{Proof :}} This proof is by induction. By definition of the value functions, we have $V^{S,k}_0(y) =0$ for all $y \in \mathscr{Y}^{S,k}$, and all $S \in \mathscr{B}$ and all $k \in \{1,2,\ldots,b_S\}$. Similarly, $V^{\mathscr{B}}_0(r) = 0$ for all $r \in \mathscr{Y}^{\mathscr{B}}$ as well. Hence, $\sum_{S \in \mathscr{B}} \sum_{k=1}^{b_S} V^{S,k}_0(r_{S,k}) = V^{\mathscr{B}}_0(r)$ for all $r \in \mathscr{R}$. Let $t \in \mathbb{N}_0$ and assume that $\sum_{S \in \mathscr{B}} \sum_{k=1}^{b_S}  V^{S,k}_t(r_{S,k}) = V^{\mathscr{B}}_t(r)$ for all $r \in \mathscr{Y}^{\mathscr{B}}$. Let $r \in \mathscr{Y}^{\mathscr{B}}$. Now, observe that
\begin{equation*} \begin{aligned} &\sum_{S \in \mathscr{B}} \sum_{k=1}^{b_S} V_{t+1}^{S,k}(r_{S,k}) \hspace{125mm} \\
=& \sum_{S \in \mathscr{B}} \sum_{k=1}^{b_S} \left(\sum_{i \in S} \left[ \lambda^*_i \min_{l \in \{0,\min\{1,r_{S,k}\}\}} \left\{ V^{S,k}_t(r_{S,k} -l) + (1-l)d_i\right\} \right. \right. \\
& \left. \left. +\mu^*_i \min_{l \in \{0,\min{\{1,C_S-r_{S,k}\}}\}} V^{S,k}_t(r_{S,k}+l) \right]+ \left( 1 - \sum_{i \in S} \left[\lambda_i^* + \mu_i^*\right] \right) V_t^{S,k}(r_{S,k}) \right) \\
& + \sum_{S \in \mathscr{B}} \sum_{k=1}^{b_S} H^{S,k}(r_{S,k})
 \end{aligned} \end{equation*}
 \begin{equation*} \begin{aligned}
=& \sum_{S \in \mathscr{B}} \sum_{k=1}^{b_S} \left(\sum_{i \in S} \left[ \lambda^*_i \min_{l \in \{0,\min\{1,r_{S,k}\}\}} \left\{ V^{S,k}_t(r_{S,k} -l) + (1-l)d_i\right\} \right. \right. \\
&  \left. \left. + \mu^*_i \min_{l \in \{0,\min\{1,C_S-r_{S,k}\}\}} V^{S,k}_t(r_{S,k}+l) \right]+ \sum_{i \in N \backslash S} \Bigg[\lambda^*_i V^{S,k}_t(r_{S,k}) + \mu^*_i V^{S,k}_t(r_{S,k}) \Bigg] \right) \\
& + \sum_{S \in \mathscr{B}} \sum_{k=1}^{b_S} H^{S,k}(r_{S,k})  \\
 =& \sum_{i \in N} \left[ \lambda^*_i \left(\sum_{S \in \mathscr{B} : i \in S} \sum_{k=1}^{b_S} \min_{l \in \{0,\min\{r_{S,k},1\}\}} \left\{V^{S,k}_t(r_{S,k} - l) + (1 - l)d_i \right\} + \sum_{S \in \mathscr{B} : i \not\in S} \sum_{k=1}^{b_S} V^{S,k}_t(r_{S,k}) \right) \right. \\
& +\mu^*_i \left( \sum_{S \in \mathscr{B} : i \in S} \sum_{k=1}^{b_S} \min_{l \in \{0,\min{\{1,C_S-r_{S,k}\}}\}} V^{S,k}_{t}(r_{S,k}+l)  + \sum_{S \in \mathscr{B} : i \not\in S} \sum_{k=1}^{b_S} V^{S,k}_t(r_{S,k})\right) \Bigg] \\
&+\sum_{S \in \mathscr{B}} \sum_{k=1}^{b_S} H^{S,k}(r_{S,k}) \\
 =&  \sum_{i \in N} \left[\lambda^*_i \cdot \left( \min_{l \in \mathscr{A}^{\mathscr{B}}_{i,-}(r)} \left\{ \sum_{S \in \mathscr{B} : i \in S} \sum_{k=1}^{b_S} V^{S,k}_t(r_{S,k} - l_{S,k}) + \sum_{S \in \mathscr{B} : i \not\in S} \sum_{k=1}^{b_S} V^{S,k}_t(r_{S,k}) +( \alpha - \vert\vert l \vert \vert_{_1})d_i \right\}  \right) \right. \\
& + \left. \mu^*_i \cdot \left( \min_{l \in \mathscr{A}^{\mathscr{B}}_{i,+}(r)} \left\{ \sum_{S \in \mathscr{B}:i \in S} \sum_{i=1}^{b_S} V^{S,k}_{t}(r_{S,k}+l_{S,k}) + \sum_{S \in \mathscr{B}:i \not\in S} \sum_{i=1}^{b_S} V^{S,k}_{t}(r_{S,k}) \right\}      \right)\right] \\
& + \sum_{S \in \mathscr{B}} \sum_{k=1}^{b_S} H^{S,k}(r_{S,k}) \\
=&  \sum_{i \in N} \left[\lambda^*_i \cdot \left( \min_{l \in \mathscr{A}^{\mathscr{B}}_{i,-}(r)} \left\{V^{\mathscr{B}}_t(r-l) +( \alpha - \vert\vert l \vert \vert_{_1})d_i \right\}\right)   + \mu^*_i \cdot \left( \min_{l \in \mathscr{A}^{\mathscr{B}}_{i,+}(r)} \left\{V_t^{\mathscr{B}}(r+l)\right\} \right) \right] \\
&+ \sum_{S \in \mathscr{B}} \sum_{k=1}^{b_S} H^{S,k}(r_{S,k}) \\
=&V^{\mathscr{B}}_{t+1}(r).
\end{aligned} \end{equation*}

\noindent The first equality holds by Lemma 1. The second equality holds as $1 - \sum_{i \in S} [\lambda_i^* + \mu_i^*]$ $=\sum_{i \in N} [ \lambda_i^* + \mu_i^*] - \sum_{i \in S} [\lambda_i^* + \mu_i^*] = \sum_{i \in N \backslash S} [\lambda_i^* + \mu_i^*]$. The third equality holds by writing $\lambda_i^*$ and $\mu^*_i$ in front of the summations. The fourth equality holds by using the definition of the $L^1$ norm, the fact that the sum of minima can be rewritten as one minimum and $\mathscr{A}_{i,+}^{\mathscr{B}}$ and  $\mathscr{A}_{i,-}^{\mathscr{B}}$ are defined such that the decisions made for all minima fit. In the fifth equality, we first use that $\sum_{S \in \mathscr{B} : i \not\in S} \sum_{k=1}^{b_S} V^{S,k}_t(r_{S,k}) =  \sum_{S \in \mathscr{B} : i \not\in S} \sum_{k=1}^{b_S} V^{S,k}_t(r_{S,k} - l_{S,k}) = \sum_{S \in \mathscr{B} : i \not\in S} \sum_{k=1}^{b_S} V^{S,k}_t(r_{S,k} + l_{S,k})$, because $l_{S,k} =0$ if $i \not \in S$, and subsequently use the induction hypothesis twice to conclude that $\sum_{S \in \mathscr{B} : i \not\in S} \sum_{k=1}^{b_S} V^{S,k}_t(r_{S,k}-l_{S,k}) +  \sum_{S \in \mathscr{B} : i \in S} \sum_{k=1}^{b_S} V^{S,k}_t(r_{S,k} - l_{S,k}) = V_t^{\mathscr{B}}(r-l)$ and $\sum_{S \in \mathscr{B} : i \not\in S} \sum_{k=1}^{b_S} V^{S,k}_t(r_{S,k}+l_{S,k}) +  \sum_{S \in \mathscr{B} : i \in S} \sum_{k=1}^{b_S} V^{S,k}_t(r_{S,k} + l_{S,k}) = V_t^{\mathscr{B}}(r+l)$. The last equality holds by Definition \ref{tweedeval}. By the principle of mathematical induction, this completes the proof.$\hfill \square$ \bigskip

\noindent \underline{Proof of Lemma 8} \bigskip

\noindent \textbf{\emph{Proof} :} This proof is by induction. By definition of the value functions $V^{\mathscr{B}}_0(r) = \widehat{V}^{\mathscr{B}}_0(r) = 0$ for all $r \in \mathscr{Y}^{\mathscr{B}}$. Let $t \in \mathbb{N}_0$ and assume that $V^{\mathscr{B}}_t(r) \ge \widehat{V}^{\mathscr{B}}_t(r)$ for all $r \in \mathscr{Y}^{\mathscr{B}}$. Let $r \in \mathscr{Y}^{\mathscr{B}}$. First, we will show that
\begin{equation} \begin{aligned}\sum_{S \in \mathscr{B}} \sum_{k=1}^{b_S} H^{S,k}(r_{S,k}) \geq H^{N}_{\alpha}(\vert \vert r \vert \vert_{_1}). \label{hafschatten} \end{aligned}
\end{equation}

\noindent For all $S \subseteq N$ and all $k \in \{1,2,\ldots,b_S\}$ we define 
 $$x_i^{S,k} = \left\{\begin{matrix} \max\{0,\min\{C_i,r_{S,k} - \sum_{j \in S: j<i}C_j\}\} & & \mbox{ for all } i \in S \\ 0 & & \mbox{ for all }i \not \in S.\end{matrix} \right.$$

\noindent Observe that $0 \leq x_i^{S,k} \leq C_i$ for all $i \in S$, all $S \subseteq N$ and all $k \in \{0,1,\ldots,b_S\}$, $\sum_{i \in N} x_i^{S,k}= r_{S,k}$ and $H^{S,k}(r_{{S,k}}) = \sum_{i \in N}h_i^* \cdot x_{i}^{S,k}$. Next, we define
$$\overline{x}_i = \sum_{S \in \mathscr{B}} \sum_{k=1}^{b_S} x_{i}^{S,k} \mbox{ for all } i \in N.$$

\noindent Moreover, we define
$$X(r) = \left\{(x_i)_{i \in N} \hspace{1mm} \vert \hspace{1mm} 0 \leq x_i \leq \alpha \cdot C_i \mbox{ for all } i \in N, x_i \in \mathbb{N}_0 \mbox{ for all } i \in N, \sum_{i \in N} x_i = \vert \vert r \vert \vert_{_{1}} \right\}.$$

\noindent Observe that $(\overline{x}_i)_{i \in N} \in X(r)$. Let $\overline{x}^* \in X(r)$ for which $\sum_{i \in N} h_i^* \cdot \overline{x}_i^* \leq \sum_{i \in N} h_i^* \cdot x_i'$ for all $x' \in X(r)$. Note that $\overline{x}^*$ exists, since $X(r)$ is a finite set. Now, observe that
\begin{equation*} \begin{aligned}\sum_{S \in \mathscr{B}} \sum_{k=1}^{b_S} H^{S,k}(r_{S,k}) &= \sum_{S \in \mathscr{B}} \sum_{k=1}^{b_S} \sum_{i \in N} h_i^* \cdot x_{i}^{S,k} \\
&= \sum_{i \in N} h_i^* \sum_{S \in \mathscr{B}} \sum_{k=1}^{b_S}  x_i^{S,k} = \sum_{i \in N} h_i^* \cdot \overline{x}_i \geq \sum_{i \in N} h_i^* \cdot \overline{x}_i^* = H^{N}_{\alpha}(\vert \vert r \vert \vert_{_1}),\end{aligned}
\end{equation*}
 
\noindent where the first, third, and last equality hold by definition. The second equality holds by changing the order of the summations. The inequality holds by definition of $\overline{x}^*$. \bigskip

\noindent This establishes (\ref{hafschatten}). Next, observe that
 
\begin{equation*} \begin{aligned} V^{\mathscr{B}}_{t+1}(r) =  &\sum_{i \in N} \left[\lambda^*_i \min_{l \in \mathscr{A}^{\mathscr{B}}_{i,-}(r)}\bigg\{ (\alpha - \vert \vert l \vert \vert_{_1}) d_i + V^{\mathscr{B}}_t(r-l) \bigg\} \right. \\
&\hspace{10mm}+ \left. \mu^*_i \min_{l \in \mathscr{A}^{\mathscr{B}}_{i,+}(r)} \left\{ V^{\mathscr{B}}_t(r+l)  \right\} \right] + \sum_{S \in \mathscr{B}} \sum_{k=1}^{b_S} H^{S,k}(r_{S,k}) \\
\geq  &\sum_{i \in N} \left[\lambda^*_i \min_{l \in \widehat{\mathscr{A}}^{\hspace{1.3mm} \mathscr{B}}_{i,-}(r)}\bigg\{ (\alpha - \vert \vert l \vert \vert_{_1}) d_i + V^{\mathscr{B}}_t(r-l) \bigg\} \right. \\
&\hspace{10mm}+ \left. \mu^*_i \min_{l \in \widehat{\mathscr{A}}^{\hspace{1.3mm} \mathscr{B}}_{i,+}(r)} \left\{ V^{\mathscr{B}}_t(r+l)  \right\} \right] + H^{N}_{\alpha}(\vert \vert r \vert \vert_{_1}) \\
\geq  &\sum_{i \in N} \left[\lambda^*_i \min_{l \in \widehat{\mathscr{A}}^{\hspace{1.3mm} \mathscr{B}}_{i,-}(r)}\bigg\{ (\alpha - \vert \vert l \vert \vert_{_1}) d_i + \widehat{V}^{\mathscr{B}}_t(r-l) \bigg\} \right. \\
&\hspace{10mm}+ \left. \mu^*_i \min_{l \in \widehat{\mathscr{A}}^{\hspace{1.3mm} \mathscr{B}}_{i,+}(r)} \left\{ \widehat{V}^{\mathscr{B}}_t(r+l)  \right\} \right] + H^{N}_{\alpha}(\vert \vert r \vert \vert_{_1}) \\
=& \widehat{V}_{t+1}^{\mathscr{B}}(r).
\end{aligned} 
\end{equation*}
 
 \noindent The first and last equality hold by definition. The first inequality holds by Lemma \ref{kakieplant} and (\ref{hafschatten}). The second inequality holds by induction hypothesis. This completes the proof by induction. $\hfill \square$ \bigskip 
  
%

\noindent  \underline{Proof of Lemma \ref{ano}} \bigskip

\noindent \textbf{\emph{Proof :}} This proof is by induction. By definition of the value functions $\widehat{V}^{\mathscr{B}}_0(r) =  V^{\alpha}_0(\vert \vert r \vert \vert_{_1}) = 0$ for all $r \in \mathscr{Y}^{\mathscr{B}}$. Let $t \in \mathbb{N}_0$ and assume that $\widehat{V}^{\mathscr{B}}_t(r) = V^{\alpha}_t(\vert \vert r \vert \vert_{_1})$ for all $r \in \mathscr{Y}^{\mathscr{B}}$. Let $r \in \mathscr{Y}^{\mathscr{B}}$. Now, observe that
\begin{equation*} \begin{aligned} \widehat{V}^{\mathscr{B}}_{t+1}(r) =&  \sum_{i \in N} \lambda^*_i \min_{l \in \widehat{\mathscr{A}}^{\hspace{1.3mm} \mathscr{B}}_{i,-}(r)}\bigg\{ (\alpha - \vert \vert l \vert \vert_1) d_i + \widehat{V}_t^{\mathscr{B}}(r-l) \bigg\} +\sum_{i \in N} \mu^*_i \min_{l \in \widehat{\mathscr{A}}^{\hspace{1mm} \mathscr{B}}_{i,+}(r)} \left\{ \widehat{V}_t^{\mathscr{B}}(r+l) \right\} \\
& + H^{N}_{\alpha}( \vert \vert r \vert \vert{_1})\\
  =&  \sum_{i \in N} \lambda^*_i \min_{z \in \{0,1,\ldots, \min\{\alpha, \vert \vert r \vert \vert_1\}\}} \left\{\min_{l \in \widehat{\mathscr{A}}^{\hspace{1.3mm} \mathscr{B}}_{i,-}(r) : \vert \vert l \vert \vert_1 = z} \bigg\{ (\alpha - z) d_i + \widehat{V}_t^{\mathscr{B}}(r-l) \bigg\} \right\} \\
  & +\sum_{i \in N} \mu^*_i \min_{ z \in \{0,1,\ldots, \min\{\alpha, \alpha \cdot C_N - \vert \vert r \vert \vert_1\}\}} \left\{\min_{l \in \widehat{\mathscr{A}}^{\hspace{1.3mm} \mathscr{B}}_{i,+}(r) : \vert \vert l \vert \vert_1 = z}  \widehat{V}_t^{\mathscr{B}}(r+l) \right\} +  H^{N}_{\alpha}( \vert \vert r \vert \vert{_1})\\
   =&   \sum_{i \in N} \lambda^*_i \min_{z \in \{0,1,\ldots, \min\{\alpha, \vert \vert r \vert \vert_1\}\}} \left\{ (\alpha - z) d_i + V^{\alpha}_t(\vert \vert r \vert \vert_{_{1}} -z) \right\}  \\
  &+ \sum_{i \in N} \mu^*_i \min_{ z \in \{0,1,\ldots, \min\{\alpha, \alpha \cdot C_N - \vert \vert r \vert \vert_1\}\}} \bigg\{  V^{\alpha}_t(\vert \vert r \vert \vert_1 +z) \bigg\} + H^{N}_{\alpha}( \vert \vert r \vert \vert{_1})\\
  =&  V^{\alpha}_{t+1}(\vert \vert r \vert \vert_{_1}).
\end{aligned} \end{equation*}

\noindent The first equality holds by Definition \ref{valuefunctionsint}. The second equality holds by rewriting the minimum as a two-step minimization. The third equality holds by the induction hypothesis. The last equality holds by Definition \ref{valuesint2}. By the principle of mathematical induction, this completes the proof. $\hfill \square$ \bigskip 

\noindent \underline{Proof of Lemma \ref{allesofniks}} \bigskip

\noindent \textbf{\emph{Proof}} : First, the value function will be rewritten. For all $j \in \{0,1,\ldots,\alpha \cdot C_N \}$ and all $t \in \mathbb{N}_0$, we have
 \begin{equation*} \begin{aligned}
V^{\alpha}_{t+1}(j) =&  \sum_{i \in N} \lambda^*_i \min_{l \in \{0,1,\ldots, \min\{\alpha,j\}\}}\bigg\{ (\alpha - l) d_i + V^{\alpha}_t(j-l) \bigg\} \\
 & +\sum_{i \in N} \mu^*_i  \min_{l \in \left\{0,\ldots,\min\left\{\alpha, \alpha \cdot C_N-j\right\}\right\}} \left\{V^{\alpha}_t(j+l) \right\}+ H^{N}_{\alpha}(j) \\
 =& \sum_{i \in N} \lambda^*_i \left[ \min_{l \in \{0,1,\ldots, \min\{\alpha,j\}\}}\bigg\{ (j-l) d_i + V^{\alpha}_t(j-l) \bigg\} + (\alpha - j) d_i \right] \\
 & +\sum_{i \in N} \mu^*_i  \min_{l \in \left\{0,\ldots,\min\left\{\alpha, \alpha \cdot C_N-j\right\}\right\}} \left\{V^{\alpha}_t(j+l) \right\}+ H^{N}_{\alpha}(j) \\
 =&  \sum_{i \in N} \lambda^*_i \left[ \min_{l \in \{\max\{0,j-\alpha\},\ldots,j\}}\bigg\{ l d_i + V^{\alpha}_t(l) \bigg\} + (\alpha - j) d_i \right] \\
 & +\sum_{i \in N} \mu^*_i  \min_{l \in \left\{j,\ldots,\min\left\{j+\alpha,\alpha \cdot C_N\right\}\right\}} \left\{V^{\alpha}_t(l) \right\}+ H^{N}_{\alpha}(j), \\
  \end{aligned} \end{equation*}

\noindent where the second equality holds as $(a-l)d_i = (\alpha -j)d_i + (j-l)d_i$ and $(\alpha-j)d_i$ is a constant. The last equality holds by substituting $j-l$ into a new variable and by substituting $j+l$ into a new variable. \bigskip

\noindent In addition, we define for all $j \in \{0,1,\ldots, \alpha \cdot C_N\}$ and all $t \in \mathbb{N}_0$

\begin{equation*} \begin{aligned} V^{\alpha_{_1}}_{t+1}(j) &= \sum_{i \in N} \lambda^*_i  \min_{l \in \{\max\{0,j-\alpha\},\ldots,j\}}\bigg\{ l d_i + V^{\alpha}_t(l) \bigg\}  \\
V^{\alpha_{_2}}_{t+1}(j) &= \sum_{i \in N} \mu^*_i  \min_{l \in \left\{j,\ldots,\min\left\{j+\alpha,\alpha \cdot C_N\right\}\right\}}\left\{V^{\alpha}_t(l) \right\}\\
V^{\alpha_{_3}}_{t+1}(j) &=  \sum_{i \in N} \lambda_i^* (\alpha - j) d_i \\
V^{\alpha_{_4}}_{t+1}(j) &= H^{N}_{\alpha}(j). \\
\end{aligned}\end{equation*}

\noindent Note that $V^{\alpha}_t(j) = V^{\alpha_{_1}}_t(j) + V^{\alpha_{_2}}_t(j)  + V^{\alpha_{_3}}_t(j)+ V^{\alpha_{_4}}_t(j)$ for all $j \in \{0,1,\ldots,\alpha \cdot C_N\}$ and all $t \in \mathbb{N}_0$. \bigskip

\noindent $(i)$ Now, we will prove by induction that $V^{\alpha}_t(j) + V^{\alpha}_t(j+2) \geq 2 \cdot V^{\alpha}_t(j+1)$ for all $j \in \{0,1,\ldots,\alpha \cdot C_N - 2\}$ and all $t \in \mathbb{N}_0$. By definition of the value functions $V^{\alpha}_0(j) = 0$ for all $j \in \{0,1,\ldots,\alpha \cdot C_N\}$. Hence, $V^{\alpha}_0(j) + V^{\alpha}_0(j+2) \geq 2 \cdot V^{\alpha}_0(j+1)$ for all $j \in \{0,1,\ldots,\alpha \cdot C_N - 2\}$. Let $t \in \mathbb{N}_0$ and assume that $V^{\alpha}_t(j) + V^{\alpha}_t(j+2) \geq 2 \cdot V^{\alpha}_t(j+1)$ for all $j \in \{0,1,\ldots,\alpha \cdot C_N - 2\}$. \bigskip

\noindent We first focus on $V_{t}^{\alpha_1}(\cdot)$, thereafter we focus on $V_t^{\alpha_2}(\cdot)$, then we focus on $V_t^{\alpha_3}(\cdot)$, and finally we focus on $V_t^{\alpha_4}(\cdot)$. We split the analysis of $V_t^{\alpha_1}(\cdot)$ in two cases. 

\noindent Let $j \in \{0,1, \ldots, \alpha-1\}$. Now, observe that
\begin{equation*} \begin{aligned}
&V^{\alpha_{_1}}_{t+1}(j) + V^{\alpha_{_1}}_{t+1}(j+2) \\
&= \sum_{i \in N} \lambda^*_i \min_{l \in \{\max\{0,j-\alpha\},\ldots,j\}}\bigg\{ l d_i + V^{\alpha}_t(l) \bigg\} + \sum_{i \in N} \lambda^*_i \min_{l \in \{\max\{0,j+2-\alpha\},\ldots,j+2\}}\bigg\{ l d_i + V^{\alpha}_t(l) \bigg\} \\
&\geq \sum_{i \in N} \lambda^*_i \min_{\substack{l_1 \in \{0,\ldots,j\} \\
 l_2 \in \{0,\ldots,j+2\}}}\bigg\{ (l_1+l_2) d_i + V^{\alpha}_t(l_1) + V^{\alpha}_t(l_2) \bigg\} \\
 &\geq 2\sum_{i \in N} \lambda^*_i \min_{l_3 \in \{0,1,\ldots,j+1\}} \{  V^{\alpha}_t(l_3) + l_3d_i \}\bigg\} \\
 &= 2\sum_{i \in N} \lambda^*_i \min_{l_3 \in \{\max\{0,j+1-\alpha\}, \ldots, j+1\}} \{  V^{\alpha}_t(l_3) + l_3d_i \}\bigg\} \\
 &= 2 V^{\alpha_{_1}}_{t+1}(j+1).
 \end{aligned} \end{equation*}

\noindent The first inequality holds as adding a (possible) term to a set from which its minimum is selected will not increase the minimum. The second inequality holds by Lemma \ref{bewijsvoormidpointconvex} 
(this lemma is at the end of this appendix) with $f(x) = V_t^{\alpha}(x)+x\cdot d_i$ for all $x \in \{0,1,\ldots, \alpha \cdot C_N\}$ (which is convex as the sum of a convex function and a linear one is still convex), $a=0$, $b= j$, $c=0$, and $d=2$.
The last but one equality holds as $j+1-\alpha \leq 0$. \bigskip

\noindent Let $j \in \{\alpha, \alpha +1, \ldots, \alpha \cdot C_N - 2\}$. Now, observe that
\begin{equation*} \begin{aligned}
&V^{\alpha_{_1}}_{t+1}(j) + V^{\alpha_{_1}}_{t+1}(j+2) \\
&=\sum_{i \in N} \lambda^*_i \min_{l \in \{\max\{0,j-\alpha\},\ldots,j\}}\bigg\{ l d_i + V^{\alpha}_t(l) \bigg\} + \sum_{i \in N} \lambda^*_i \min_{l \in \{\max\{0,j+2-\alpha\},\ldots,j+2\}}\bigg\{ l d_i + V^{\alpha}_t(l) \bigg\} \\
&= \sum_{i \in N} \lambda^*_i \min_{\substack{l_1 \in \{j-\alpha,\ldots,j\} \\
 l_2 \in \{j+2-\alpha,\ldots,j+2\}}}\bigg\{ (l_1+l_2) d_i + V^{\alpha}_t(l_1) + V^{\alpha}_t(l_2) \bigg\} 
  \end{aligned} \end{equation*}
  \begin{equation*} \begin{aligned}
 &\geq 2\sum_{i \in N} \lambda^*_i \min_{l_3 \in \{j+1-\alpha, \ldots,j+1\}} \{  V^{\alpha}_t(l_3) + l_3d_i \}\bigg\} \hspace{74mm} \\
 &= 2\sum_{i \in N} \lambda^*_i \min_{l_3 \in \{\max\{0,j+1-\alpha\}, \ldots,j+1\}} \{  V^{\alpha}_t(l_3) + l_3d_i \}\bigg\}. \\
 &= 2 V^{\alpha_{_1}}_{t+1}(j+1).
 \end{aligned} \end{equation*}

\noindent The first inequality holds by Lemma \ref{bewijsvoormidpointconvex} with $f(x) = V_t^{\alpha}(x)+x\cdot d_i$ for all $x \in \{0,1,\ldots, \alpha \cdot C_N\}$ (which is convex as the sum of a convex function and a linear one is still convex), $a=j-\alpha$, $b=j$, $c=2$, and $d=2$. The last but one equality holds as $j+1-\alpha   \geq 0$. \bigskip

 \noindent So, for all $j \in \{0,1,\ldots,\alpha \cdot C_N - 2\}$ it holds that $V^{\alpha_{_1}}_{t+1}(j) + V^{\alpha_{_1}}_{t+1}(j+2) \geq 2 V^{\alpha_{_1}}_{t+1}(j+1)$. \bigskip

\noindent We now focus on $V^{\alpha_2}_t(\cdot)$ and split the analysis in two cases.

\noindent Let $j \in \{0,1,\ldots,\alpha \cdot (C_N-1) - 2\}$. Now, observe that
\begin{equation*} \begin{aligned}
V^{\alpha_{_2}}_{t+1}(j) + V^{\alpha_{_2}}_{t+1}(j+2) &=  \sum_{i \in N} \mu^*_i  \min_{l \in \left\{j,\ldots,\min\left\{j+\alpha,\alpha \cdot C_N\right\}\right\}}V^{\alpha}_t(l) + \sum_{i \in N} \mu^*_i  \min_{l \in \left\{j+2,\ldots,\min\left\{j+2+\alpha,\alpha \cdot C_N\right\}\right\}}V^{\alpha}_t(l)\\
&= \sum_{i \in N} \mu^*_i \min_{\substack{l_1 \in \{j,\ldots,j+\alpha \} \\
 l_2 \in \{j+2,\ldots, j+\alpha+ 2\}}} V^{\alpha}_t(l_1) + V^{\alpha}_t(l_2) \\
&\geq 2 \sum_{i \in N} \mu^*_i \min_{l_3 \in \left\{j+1,\ldots, j+\alpha+1\right\}} V^{\alpha}_t(l_3)  \\
&= 2 \sum_{i \in N} \mu^*_i \min_{l_3 \in \left\{j+1,\ldots, \min\left\{j+\alpha+1, \alpha \cdot C_N\right\}\right\}} V^{\alpha}_t(l_3) \\
& = 2 \cdot V^{\alpha_2}_{t+1}(j+1)
\end{aligned} \end{equation*}

\noindent The inequality holds by Lemma \ref{bewijsvoormidpointconvex} with $f(x) = V^{\alpha}_t(x)$ for all $x \in \{0,1,\ldots, \alpha \cdot C_N\}$, $a= j$, $b=j+\alpha$, $c = d= 2$. The last but one equality holds as $j \leq \alpha \cdot (C_N-1)-2$. \bigskip

\noindent Let $j \in \{\alpha \cdot (C_N-1) - 1, \alpha \cdot (C_N-1), \ldots, \alpha \cdot C_N -2\}$. Now, observe that
\begin{equation*} \begin{aligned}
V^{\alpha_{_2}}_{t+1}(j) + V^{\alpha_{_2}}_{t+1}(j+2) &= \sum_{i \in N} \mu^*_i  \min_{l \in \left\{j,\ldots,\min\left\{j+\alpha,\alpha \cdot C_N\right\}\right\}}V^{\alpha}_t(l) + \sum_{i \in N} \mu^*_i  \min_{l \in \left\{j+2,\ldots,\min\left\{j+2+\alpha,\alpha \cdot C_N\right\}\right\}}V^{\alpha}_t(l) \\
&\geq \sum_{i \in N} \mu^*_i \min_{\substack{l_1 \in \{j,\ldots,\alpha \cdot C_N\} \\
 l_2 \in \{j+2,\ldots,\alpha \cdot C_N\}}} \left\{V^{\alpha}_t(l_1) + V^{\alpha}_t(l_2) \right\} \\
&\geq 2 \sum_{i \in N} \mu^*_i \min_{l_3 \in \left\{j+1,\ldots, \alpha \cdot C_N\right\}} V^{\alpha}_t(l_3)  \\
&= 2 \sum_{i \in N} \mu^*_i \min_{l_3 \in \left\{j+1,\ldots, \min\left\{j+1+\alpha, \alpha \cdot C_N\right\}\right\}} V^{\alpha}_t(l_3) \\
& = 2 \cdot V^{\alpha_2}_{t+1}(j+1)
\end{aligned} \end{equation*}

\noindent The first inequality holds as adding a possible term to a set from which its minimum is selected will not increase the minimum. The second inequality holds by Lemma \ref{bewijsvoormidpointconvex} with $f(x) = V^{\alpha}_t(x)$ for all $x \in \{0,1,\ldots, \alpha \cdot C_N\}$, $a=j$, $b= \alpha \cdot C_N$, $c=2$, and $d=0$. The last but one equality holds as $j+1+\alpha \geq \alpha \cdot C_N$.

\bigskip

 \noindent So, for all $j \in \{0,1,\ldots,\alpha \cdot C_N - 2\}$ it holds that $V^{\alpha_{_2}}_{t+1}(j) + V^{\alpha_{_2}}_{t+1}(j+2) \geq 2 V^{\alpha_{_2}}_{t+1}(j+1)$. \bigskip

\noindent We now focus on $V^{\alpha_3}_t(\cdot)$.

\noindent Let $j \in \{0,1,\ldots,\alpha \cdot C_N - 2\}$. Now, observe that
\begin{equation*} \begin{aligned}
V^{\alpha_{_3}}_{t+1}(j) + V^{\alpha_{_3}}_{t+1}(j+2) &=  \sum_{i \in N} \lambda^*_i(\alpha-j)d_i + \sum_{i \in N} \lambda^*_i(\alpha-(j+2))d_i \\
&=  2\sum_{i \in N} \lambda_i^* (\alpha - (j+1))d_i \\ 
&= 2 \cdot V^{\alpha_{_3}}_{t+1}(j+1).
\end{aligned} \end{equation*}

\noindent So, for all $j \in \{0,1,\ldots,\alpha \cdot C_N - 2\}$ it holds that $V^{\alpha_{_3}}_{t+1}(j) + V^{\alpha_{_3}}_{t+1}(j+2) \geq 2 V^{\alpha_{_3}}_{t+1}(j+1)$. \bigskip

Finally, we focus on $V^{\alpha_4}_t(\cdot)$.

\noindent Let $j \in \{0,1,\ldots, \sum_{i \in N} \alpha \cdot C_i - 1\}$. Then, using Definition \ref{def:holdingcosts}, we have $V_{t+1}^{\alpha_4}(j+1) - V_{t+1}^{\alpha_4}(j) = h_{k_j}^*$ for some $k_j \in N$. Because spare parts are stored in the cheapest location possible (i.e., first in location 1 until  capacity is reached, then location 2 until capacity is reached, and so on), we have $h_{k_{j'}}^* \leq h_{k_{{j'}+1}}^*$ for all $j' \in \{0,1,\ldots, \sum_{i \in N} \alpha \cdot C_i - 2\}$. 
Consequently, for a given $j \in \{0,1,\ldots, \sum_{i \in N} \alpha \cdot C_i - 2\}$, we obtain
\begin{equation*} \begin{aligned} 2 \cdot V^{\alpha_4}_{t+1}(j+1) - \left( V^{\alpha_4}_{t+1}(j) + V^{\alpha_4}_{t+1}(j+2)\right) &=
\left( V^{\alpha_4}_{t+1}(j+1) - V^{\alpha_4}_{t+1}(j)\right) - \left( V^{\alpha_4}_{t+1}(j+2) - V^{\alpha_4}_{t+1}(j+1)\right) \\
&= h_{k_j}^* - h_{k_{j+1}}^* \leq 0,\end{aligned} \end{equation*}
where the inequality holds, because $h_{k_j}^* \leq h_{k_{j+1}}^*$.

\noindent Hence, for all $j \in \{0,1,\ldots, \alpha \cdot C_N-2\}$ it holds that $V^{\alpha_4}_{t+1}(j) + V^{\alpha_4}_{t+1}(j+2) \geq 2 \cdot V^{\alpha}_{t+1}(j+1)$. \bigskip

\noindent We conclude that, for all $j \in \{0,1,\ldots, \alpha \cdot C_N - 2\}$, we have
\begin{equation*} \begin{aligned}
V^{\alpha}_{t+1}(j) + V^{\alpha}_{t+1}(j+2) =&  V^{\alpha_{_1}}_{t+1}(j) + V^{\alpha_{_2}}_{t+1}(j) +  V^{\alpha_{_3}}_{t+1}(j) +  V^{\alpha_{_4}}_{t+1}(j) \\
 &+ V^{\alpha_{_1}}_{t+1}(j+2) + V^{\alpha_{_2}}_{t+1}(j+2) + V^{\alpha_{_3}}_{t+1}(j+2)+  V^{\alpha_{_4}}_{t+1}(j+2)\\
\geq& 2 \cdot V^{\alpha_{_1}}_{t+1}(j+1) + 2 V^{\alpha_{_2}}_{t+1}(j+1) + 2 V^{\alpha_{_3}}_{t+1}(j+1) + 2 V^{\alpha_{_4}}_{t+1}(j+1)\\
=& 2 \cdot V^{\alpha}_{t+1}(j+1).  \end{aligned} \end{equation*}

  \noindent $(ii)$ \noindent Next, we will prove by induction that  $V^{\alpha}_t(k + j) + V^{\alpha}_t(k  + j+2) = 2\cdot V^{\alpha}_t(k + j +1)$ for all $j \in \{0,1,\ldots,\alpha -2\}$ and all $k \in \{0,\alpha, 2\alpha,\ldots,(C_N-1)\alpha\}$. In doing so, we first observe that by linearity we can write $V^{\alpha_1}_{t}(j) = \sum_{i \in N} V^{\alpha_1,i}_{t}(j)$
with $V^{\alpha_1,i}_{t+1}(j) = \lambda^*_i \bigg[\min_{l \in \{\max\{0,j - \alpha\},\ldots,j\}}\bigg\{ l d_i + V^{\alpha}_{t}(l) \bigg\}  \bigg]$ \noindent for all $i \in N$, all $j \in \{0,1,\ldots \alpha \cdot C_N\}$ and all $t \in \mathbb{N}_{0}$.

  By definition of the value functions $V^{\alpha}_0(j) = 0$ for all $j \in \{0,1,\ldots,\alpha \cdot C_N\}$. Hence, $V^{\alpha}_0(k + j) + V^{\alpha}_0(k+j+2) = 2 \cdot V^{\alpha}_0(k+j+1)$ for all $j \in \{0,1,\ldots,\alpha -2\}$ and all $k \in \{0,\alpha,2 \alpha,\ldots,(C_N-1)\alpha\}$. Let $t \in \mathbb{N}_0$ and assume that $V^{\alpha}_t(k+j) + V^{\alpha}_t(k+j+2) = 2 \cdot V^{\alpha}_t(k+j+1)$ for all $j \in \{0,1,\ldots,\alpha-2\}$ and all $k \in \{0,\alpha,2\alpha,\ldots, (C_N-1) \alpha\}$. \bigskip

We first focus on $V^{\alpha_1}_t(\cdot)$. 
Observe that the function described by $V^{\alpha}_t(j) +j\cdot d_i$ for all $j \in \{0,1,\ldots, \alpha \cdot C_N\}$ is convex as the first term is convex by $(i)$, the second term is linear, and the combination of a convex and linear term remains convex. 
 By our induction hypothesis, it holds that $V^{\alpha}_t(k+j) + V^{\alpha}_t(k+j+2) = 2 \cdot V^{\alpha}_t(k+j+1)$ for all $j \in \{0,1,\ldots, \alpha -2\}$ and all $k \in \{0,\alpha,2 \alpha\ldots,   (C_N-1) \alpha\}$, which implies that the function described by $V^{\alpha}_t(j) + j \cdot d_i$ for all $j \in \{0,1,\ldots, \alpha \cdot C_N\}$ is also piecewise linear. So, there exists a $p_{i,t} \in \{0,\alpha,2\alpha,\ldots,\ldots,C_N\alpha\}$ for which it holds that $V_t^{\alpha}(p_{i,t}) + p_{i,t}\cdot d_i \leq V_t^{\alpha}(j) + j \cdot d_i$ for all $j \in \{0,1,\ldots,\alpha \cdot C_N\}$ and all $i \in N$.\footnote{We would like to explicitly highlight the dependency of $p$ on $i$ and $t$, because it was not made explicit in the work of \cite{schlicherMOR}. Their proof is correct, taking this silently assumed dependency into account.} Let $i \in N$, 
fix such a $p_{i,t}$ and denote it by $p_{i,t}^* \in \{0,\alpha ,\ldots, C_N \cdot \alpha\}$. Then, for all $k \in \{0,\alpha ,\ldots,p_{i,t}^*-\alpha\}$ and all $j \in \{0,1,\ldots,\alpha\}$ it holds that
\begin{equation} \begin{aligned}  \label{laatsteverg1}  V^{\alpha_1,i}_{t+1}(k+j) 
  =
  \lambda^*_i \bigg[\min_{l \in \{\max\{0,k+j-\alpha\},\ldots,k+j\}}\bigg\{ l d_i + V^{\alpha}_t(l) \bigg\}  \bigg] 
  =   
  \lambda^*_i \bigg[ (k+j)d_i +  V^{\alpha}_t(k+j)  \bigg],\\
   \end{aligned} \end{equation}
\noindent where the second
equality holds, because the minimum of the convex function described by $V^{\alpha}_t(j')+ j' \cdot d_i$ for all $j' \in \{0,1,\ldots, \alpha \cdot C_N\}$ is attained at $j'=p_{i,t}^*$ and so the minimum after the first equality is attained at $l=k+j$. \bigskip

\noindent For $k = p_{i,t}^*$ and all $j \in \{0,1,\ldots,\alpha\}$ it holds that
\begin{equation} \begin{aligned} V^{\alpha_1,i}_{t+1}(k+j) 
  =
  \lambda^*_i \bigg[\min_{l \in \{\max\{0,k+j-\alpha\},\ldots,k+j\}}\bigg\{ l d_i + V^{\alpha}_t(l) \bigg\} \bigg] 
=  
\lambda^*_i \bigg[ k \cdot d_i + V^{\alpha}_t(k) \bigg],
   \end{aligned} \label{laatsteverg2} \end{equation}
\noindent where the second
equality results from the fact that the minimum of the convex function described by $V^{\alpha}_t(j')+ j' \cdot d_i$ for all $j' \in \{0,1,\ldots, \alpha \cdot C_N\}$ is attained at $j' =p_{i,t}^*$ and so the minimum after the first equality is attained at $l=k
$. \bigskip

\noindent For all $k \in \{p_{i,t}^*+\alpha,p_{i,t}^*+2 \cdot \alpha,\ldots,(C_N-1) \cdot \alpha \}$ and all $j \in \{0,1,\ldots,\alpha\}$ it holds that
  \begin{equation}
  \begin{aligned} 
  V^{\alpha_1,i}_{t+1}(k+j) &
  =
  \lambda^*_i \bigg[\min_{l \in \{\max\{0,k+j-\alpha\},\ldots,k+j\}}\bigg\{ l d_i + V^{\alpha}_t(l) \bigg\} \bigg] \\
  & 
  =  
  \lambda^*_i \bigg[ (k+j-\alpha)d_i + V^{\alpha}_t(k+j-\alpha) \bigg], \end{aligned}  \label{laatsteverg3} \end{equation}

\noindent where the second
equality results from the fact that the minimum of the convex function described by $V^{\alpha}_t(j')+ j' \cdot d_i$ for all $j' \in \{0,1,\ldots, \alpha \cdot C_N\}$ is attained at $j'=p_{i,t}^*$ and
so the minimum after the first equality in (\ref{laatsteverg3}) is attained at $l=k+j-\alpha$. \bigskip

\noindent So, for all $k \in \{0,1,\ldots,p_{i,t}^* - \alpha\}$ and all $j \in \{0,1,\ldots,\alpha-2\}$, we have
\begin{equation} \begin{aligned}
\label{ditevenproberen1}
&V^{\alpha_1,i}_{t+1}(k+j) + V^{\alpha_1,i}_{t+1}(k+j+2)  \\
 =& 
 \lambda^*_i \bigg[ V^{\alpha}_t(k+j) + (k+j) \cdot d_i \bigg] + 
 \lambda^*_i \bigg[ V^{\alpha}_t(k+j+2) + (k+j+2) \cdot d_i \bigg]  \\
     =&
     2 
     \lambda^*_i \bigg[ V^{\alpha}_t(k+j+1) + (k+j+1) \cdot d_i \bigg] \\
     =& 2 V^{\alpha_1,i}_{t+1}(k+j+1).
\end{aligned}, \end{equation}
\noindent where the first and last 
 equality hold by (\ref{laatsteverg1}) and the second one by the induction hypothesis.\bigskip

\noindent For $k =p_{i,t}^*$ and all $j \in \{0,1,\ldots,\alpha-2\}$ it holds that
\begin{equation} \begin{aligned}
\label{ditevenproberen2}
&V^{\alpha_1,i}_{t+1}(k+j) + V^{\alpha,i}_{t+1}(k+j+2) \\
 =
 &
 \lambda^*_i \bigg[ V^{\alpha_1}_t(k) + k \cdot d_i \bigg] + \sum_{i \in N} \lambda^*_i \bigg[ V^{\alpha}_t(k) + k \cdot d_i \bigg]  \\
     =& 2 
     \lambda^*_i \bigg[ V^{\alpha}_t(k) + k \cdot d_i \bigg] \\
     =& 2 V^{\alpha_1,i}_{t+1}(k+j+1),
\end{aligned} \end{equation}
\noindent where the first and last equality hold by (\ref{laatsteverg2}). The second equality holds by the induction hypothesis.\bigskip

\noindent For all $k \in \{p^*_{i,t}+\alpha, p^*_{i,t}+2 \cdot \alpha,\ldots, (C_N - 1)\cdot \alpha\}$ and all $j \in \{0,1,\ldots,\alpha-2\}$ it holds that
\begin{equation} \begin{aligned} \label{ditevenproberen3}
&V^{\alpha_1,i}_{t+1}(k+j) + V^{\alpha_1,i}_{t+1}(k+j+2) \\
 =
 &
 \lambda^*_i \bigg[ V^{\alpha}_t(k+j -\alpha) +(k+j-\alpha)d_i \bigg] + 
 \lambda^*_i \bigg[ V^{\alpha}_t(k+j+2-\alpha) +(k+j-\alpha+2)d_i \bigg]  \\
     =& 2 
     \lambda^*_i \bigg[ V^{\alpha}_t(k+j+1-\alpha) +(k+j-\alpha+1)d_i\bigg] \\
     =& 2 V^{\alpha_1,i}_{t+1}(k+j+1),
\end{aligned} \end{equation}
\noindent where the first and last equality hold by (\ref{laatsteverg3}) and the second equality by the induction hypothesis. \bigskip

\noindent By exploiting (\ref{ditevenproberen1}), (\ref{ditevenproberen2}), and (\ref{ditevenproberen3}) for all $i \in N$, we can conclude that for all $k \in \{0,\alpha,\ldots,\alpha\cdot (C_N-1)\}$ and all $j \in \{0,1,\ldots,\alpha-2\}$ we have
\begin{equation} \label{deeleenvanheheel1} V^{\alpha_1,i}_t(k+j) + V^{\alpha_1,i}_t(k+j+2) = 2 \cdot V^{\alpha_1,i}_t(k+j+1).\end{equation}
\noindent As $V^{\alpha_1}_t(j) = \sum_{i \in N}V_t^{\alpha_i,i}(j)$ for all $j \in \{0,1,\ldots, \alpha \cdot C_N\}$, we conclude that for all $k \in \{0,\alpha,\ldots,\alpha\cdot (C_N-1)\}$ and all $j \in \{0,1,\ldots,\alpha-2\}$
\begin{equation} \label{deeleenvanheheel} V^{\alpha_1}_t(k+j) + V^{\alpha_1}_t(k+j+2) = 2 \cdot V^{\alpha_1}_t(k+j+1),\end{equation}

\noindent  We now focus on $V^{\alpha_2}_t(\cdot)$.  Recall that $V^{\alpha}_t(\cdot)$ is convex by $(i)$. By our induction hypothesis, it holds that $V^{\alpha}_t(k+j) + V^{\alpha}_t(k+j+2) = 2 \cdot V^{\alpha}_t(k+j+1)$ for all $k \in \{0,\alpha,2 \alpha\ldots,   (C_N-1) \alpha\}$ and all $j \in \{0,1,\ldots, \alpha -2\}$. So, there exists a $p_{t} \in \{0,\alpha,2\alpha,\ldots,\ldots,C_N\cdot \alpha\}$ for which it holds that $V_t^{\alpha}(p_{t})  \leq V_t^{\alpha}(j) $ for all $j \in \{0,1,\ldots,\alpha \cdot C_N\}$ and all $i \in N$. Fix such a $p_t$ and denote it by 
$\tilde{p}_t \in \{0,\alpha ,\ldots, C_N \cdot \alpha\}$. 
Now, for all $k \in \{0,\alpha ,\ldots, 
      \tilde{p}_t-2\alpha\}$ and all $j \in \{0,1,\ldots,\alpha\}$ it holds that
  \begin{equation} \begin{aligned} V^{\alpha_2}_t(k+j) & = \sum_{i \in N} \mu^*_i  \min_{l \in \left\{k+j,\ldots,\min\left\{k+j+\alpha,\alpha \cdot C_N\right\}\right\}}V^{\alpha}_t(l) \\
  &=   \sum_{i \in N} \mu^*_i  V^{\alpha}_t(k+j+\alpha),
    \label{laatsteverg11}  \end{aligned} \end{equation}

\noindent where the second equality holds as the minimum of convex function $V^{\alpha}_t(\cdot)$ is attained at $\tilde{p}_t$ and so the minimum after the first equality is attained at $l=k+j+\alpha$. \bigskip

\noindent For $k = \tilde{p}_t- \alpha$ and all $j \in \{0,1,\ldots,\alpha\}$ it holds that
  \begin{equation} \begin{aligned} V^{\alpha_2}_t(k+j) &= \sum_{i \in N} \mu^*_i  \min_{l \in \left\{k+j,\ldots,\min\left\{k+j+\alpha,\alpha \cdot C_N\right\}\right\}}V^{\alpha}_t(l) \\
&=  \sum_{i \in N} \mu^*_i V^{\alpha}_t(\tilde{p}_t), \end{aligned} \label{laatsteverg12} \end{equation}

\noindent where the second equality results from the fact that the minimum of convex function $V^{\alpha}_t(\cdot)$ is attained at $\tilde{p}_t$ and so the minimum after the first equality is attained at $l=\tilde{p}_t$. \bigskip

\noindent For all $k \in \{\tilde{p}_t,\tilde{p}_t+ \alpha,\ldots,(C_N-1) \cdot \alpha \}$ and all $j \in \{0,1,\ldots,\alpha\}$ it holds that
  \begin{equation} \label{laatsteverg13}
  \begin{aligned} V^{\alpha_2}_t(k+j) &= \sum_{i \in N} \mu^*_i  \min_{l \in \left\{k+j,\ldots,\min\left\{k+j+\alpha,\alpha \cdot C_N\right\}\right\}}V^{\alpha}_t(l) \\
  &= \sum_{i \in N} \mu^*_i  V^{\alpha}_t(k+j), \end{aligned}
  \end{equation}

\noindent where the second equality results from the fact that the minimum of convex function $V^{\alpha}_t(\cdot)$ is attained at $\tilde{p}_t$ and so the minimum after the first equality is attained at $l=k+j$. \bigskip


\noindent So, for all $k \in \{0,1,\ldots,\tilde{p}_t - 2\alpha\}$ and all $j \in \{0,1,\ldots,\alpha-2\}$, we have\begin{equation*} \begin{aligned}
&V^{\alpha_2}_{t+1}(k+j) + V^{\alpha_2}_{t+1}(k+j+2) \\
 =&\sum_{i \in N} \mu^*_i  V^{\alpha}_t(k+j+\alpha) + \sum_{i \in N} \mu^*_i  V^{\alpha}_t(k+j+2+\alpha) \\
     =& 2 \sum_{i \in N} \mu^*_i  V^{\alpha}_t(k+j+1+\alpha) \\
     =& 2 V^{\alpha_2}_{t+1}(k+j+ 1),
\end{aligned} \end{equation*}

\noindent where the first and last equality hold by (\ref{laatsteverg11}) and the second equality holds by the induction hypothesis.\bigskip

\noindent For $k =\tilde{p}_t-\alpha$ and all $j \in \{0,1,\ldots,\alpha-2\}$ it holds that
\begin{equation*} \begin{aligned}
&V^{\alpha_2}_{t+1}(k+j) + V^{\alpha_2}_{t+1}(k+j+2) \\
 =&\sum_{i \in N} \mu^*_i  V^{\alpha}_t(\tilde{p}_t) + \sum_{i \in N} \mu^*_i  V^{\alpha}_t(\tilde{p}_t) \\
     =& 2 V^{\alpha_2}_{t+1}(k+j+1),
\end{aligned} \end{equation*}

\noindent where the first and second equality hold by (\ref{laatsteverg12}). 
\bigskip

\noindent For all $k \in \{\tilde{p}_t, \tilde{p}_t+\alpha,\ldots, (C_N - 1)\cdot \alpha\}$ and all $j \in \{0,1,\ldots,\alpha-2\}$ it holds that
\begin{equation*} \begin{aligned}
&V^{\alpha_2}_{t+1}(k+j) + V^{\alpha_2}_{t+1}(k+j+2) \\
 =&\sum_{i \in N} \mu^*_i  V^{\alpha}_t(k+j) + \sum_{i \in N} \mu^*_i  V^{\alpha}_t(k+j+2) \\
     =& 2 \sum_{i \in N} \mu^*_i  V^{\alpha}_t(k+j+1) \\
     =& 2 V^{\alpha_2}_{t+1}(k+j+1),
\end{aligned} \end{equation*}

\noindent where the first and last equality hold by (\ref{laatsteverg13}) and the second equality holds by the induction hypothesis. \bigskip

\noindent We conclude that for all $k \in \{0,\alpha,\ldots,\alpha\cdot (C_N-1)\}$ and all $j \in \{0,1,\ldots,\alpha-2\}$ we have
\begin{equation} \label{deeltweevangeheel} V^{\alpha_2}_{t+1}(k+j) + V^{\alpha_2}_{t+1}(k+j+2) = 2 \cdot V^{\alpha_2}_{t+1}(k+j+1).\end{equation}

\noindent We now focus on $V^{\alpha_3}_t(\cdot)$. For all $k \in \{0,1,\ldots,(C_N-1)\alpha\}$ and all $j \in \{0,1,\ldots,\alpha-2\}$ we have
\begin{equation} \begin{aligned} \label{deeldrievangeheel} V^{\alpha_3}_{t+1}(k+j) + V^{\alpha_3}_{t+1}(k+j+2) =& 
\sum_{i \in N}\lambda_i^* (\alpha - j)d_i 
+ \sum_{i \in N}\lambda_i^* (\alpha - (j+2))d_i \\
=& 
2 \sum_{i \in N}\lambda_i^* (\alpha - (j+1))d_i\\
=& 2 V^{\alpha_3}_{t+1}(k+j+1). \end{aligned}\end{equation}

\noindent Finally, we focus on $V^{\alpha_4}_t(\cdot)$. Recall that for any 
 $k \in \{0,\alpha,\ldots,(C_N-1)\alpha\}$ there exists a $z_k \in N$ such that $V_{t+1}^{\alpha_4}(k+j+1) - V_{t+1}^{\alpha_4}(k+j) = h_{z_k}^*$ for all $j \in \{0,1,\ldots,\alpha-2\}$. Consequently, for a given $k \in \{0,\alpha,\ldots,(C_N-1)\alpha\}$ and all $j \in \{0,1,\ldots,\alpha-2\}$, we have
\begin{equation*} \begin{aligned} &2 \cdot V^{\alpha_4}_{t+1}(k+j+1) - \left( V^{\alpha_4}_{t+1}(k+j) + V^{\alpha_4}_{t+1}(k+j+2)\right) \\
&=
\left( V^{\alpha_4}_{t+1}(k+j+1) - V^{\alpha_4}_{t+1}(k+j)\right) - \left( V^{\alpha_4}_{t+1}(k+j+2) - V^{\alpha_4}_{t+1}(k+j+1)\right) \\
&= h_{z_k}^* - h_{z_k}^* = 0, \end{aligned} \end{equation*}
\noindent and thus
\begin{equation*} \begin{aligned} 
V^{\alpha_4}_{t+1}(k+j) + V^{\alpha_4}_{t+1}(k+j+2) 
= 2 \cdot V^{\alpha_4}_{t+1}(k+j+1). 
\end{aligned} \end{equation*} 


In conclusion, we have for all $k \in \{0,1,\ldots,(C_N-1)\alpha\}$ and all $j \in \{0,1,\ldots,\alpha-2\}$
\begin{equation*} \begin{aligned} V^{\alpha}_{t+1}(k+j) + V^{\alpha}_{t+1}(k+j+2) =& V^{\alpha_1}_{t+1}(k+j) + V^{\alpha_2}_{t+1}(k+j) + V^{\alpha_3}_{t+1}(k+j) + V^{\alpha_4}_{t+1}(k+j)\\
 &+ V^{\alpha_1}_{t+1}(k+j+2) + V^{\alpha_2}_{t+1}(k+j+2)\\
 & + V^{\alpha_3}_{t+1}(k+j+2) + V^{\alpha_4}_{t+1}(k+j+2) \\
=& 2 \cdot V^{\alpha_1}_{t+1}(k+j+1) + 2 \cdot V^{\alpha_2}_{t+1}(k+j+1) \\
&+ 2 \cdot V^{\alpha_3}_{t+1}(k+j+1)) + 2 \cdot V^{\alpha_4}_{t+1}(k+j+1)) \\
=& 2 \cdot V^{\alpha}_{t+1}(k+j+1),
\end{aligned} \end{equation*}

\noindent By the principle of mathematical induction, this completes the proof.   $\hfill \square$  \bigskip

\noindent \underline{Proof of Lemma \ref{alphading}} \bigskip

\noindent \textbf{\emph{Proof :}} Based on Definition \ref{valuesint2} and $(i)$ and $(ii)$ of Lemma \ref{allesofniks}, it follows directly that for all $j \in \{0,\alpha,\ldots, \alpha \cdot C_N\}$ and all $t \in \mathbb{N}_0$ it holds that
\begin{equation} \begin{aligned} \label{enerlaatstebewijslabel} V^{\alpha}_{t+1}(j) =&   \sum_{i \in N} \left[\lambda^*_i \min_{l \in \{0,\min\{j,\alpha\}\}} \left\{V^{\alpha}_t(j-l) + (\alpha-l)d_i \right\}\right]  \\
 &+ \sum_{i \in N} \left[\mu^*_i \min_{l \in \left\{0, \min\left\{ \alpha, C_N-j\right\} \right\}}\left\{V^{\alpha}_t(j+l)\right\}\right] + H^{N}_{\alpha}(j)\end{aligned} \end{equation}

\noindent Moreover, for all $j \in \{0,\alpha, \ldots, \alpha \cdot C_N\}$ it holds that
\begin{equation} \begin{aligned} H^{N}_{\alpha}(j) &= \sum_{i \in N} h_i^* \max\{0, \min\{\alpha \cdot C_i, j - \sum_{k \in N: k < i} \alpha \cdot C_k\}\} \\
&= \sum_{i \in N} h_i^* \max\{\alpha \cdot 0, \min\{\alpha \cdot C_i, \alpha \cdot \frac{j}{\alpha} - \alpha \cdot \sum_{k \in N: k < i} C_k\}\} \\
&= \alpha \sum_{i \in N} h_i^* \max\{0, \min\{ C_i,   \frac{j}{\alpha} - \sum_{k \in N: k < i} C_k\}\} \\
&=\alpha \cdot H^{N}\left(\frac{j}{\alpha}\right)
\end{aligned} \label{hvergelikingalfa} \end{equation}

\noindent The first and last equality hold by definition. In the second equality, we use that $\alpha$ is a (strictly positive) constant. The third equality holds as $\alpha$ can be taken outside the minimum and maximum.   \bigskip

  \noindent We will prove by induction that $V^{\alpha}_t(j) = \alpha \cdot V_t^N\left(\frac{j}{\alpha}\right)$ for all $j \in \{0,\alpha,\ldots,C_N \cdot \alpha\}$. By definition of the value functions $V^{\alpha}_0(j) = V_0^N(j)=0$ for all $j \in \{0,1,\ldots,C_N\}$ and so $V^{\alpha}_0(j) = \alpha \cdot V_0^N\left( \frac{j}{\alpha}\right)$ for all $j \in \{0,\alpha,\ldots,C_N \cdot \alpha\}$. Let $t \in \mathbb{N}_0$ and assume that $V^{\alpha}_t(j) = \alpha \cdot V_t^N\left(\frac{j}{\alpha}\right)$ for all $j \in \{0,\alpha,\ldots,C_N \cdot \alpha\}$. \bigskip

  \noindent Let $j \in \{0,\alpha,\ldots,C_N \cdot \alpha\}$. Then, observe that
  \begin{equation*} \begin{aligned}
  V^{\alpha}_{t+1}(j) =&  \sum_{i \in N} \left[\lambda^*_i \min_{l \in \{0,\min\{j,\alpha\}\}} \left\{V^{\alpha}_t(j-l) + (\alpha-l)d_i \right\}\right]  \\
  & + \sum_{i \in N} \left[\mu^*_i \min_{l \in \left\{0, \min\left\{ \alpha, \alpha C_N-j\right\} \right\}}\left\{V^{\alpha}_t(j+l)\right\}\right] + H^{N}_{\alpha}(j) \\
  =&  \sum_{i \in N} \left[\lambda^*_i \min_{l \in \{0,\min\{j,\alpha\}\}} \left\{\alpha \cdot V_t^N\left(\frac{j-l}{\alpha}\right) +  (\alpha-l)d_i \right\}\right]  \\
  & + \sum_{i \in N} \left[\mu^*_i \min_{l \in \left\{0, \min\left\{ \alpha, \alpha C_N-j\right\} \right\}}\left\{\alpha \cdot  V^{N}_t\left(\frac{j+l}{\alpha}\right)\right\}\right] + \alpha \cdot H^{N}\left(\frac{j}{\alpha}\right) \\
  =&  \sum_{i \in N} \left[\lambda^*_i \min_{z \in \{0,\min\{\frac{j}{\alpha},1\}\}} \left\{\alpha \cdot V_t^N\left(\frac{j}{\alpha} - z\right) +  \alpha \cdot(1-z)d_i \right\}\right]  \\
  & + \sum_{i \in N} \left[\mu^*_i \min_{z \in \left\{0, \min\left\{1, C_N - \frac{j}{\alpha}\right\} \right\}}\left\{\alpha \cdot  V^{N}_t\left(\frac{j}{\alpha} + z \right)\right\}\right] + \alpha \cdot H^{N}\left(\frac{j}{\alpha}\right) \\
      =& \alpha \cdot \left(  \sum_{i \in N} \left[\lambda^*_i \min_{z \in \{0,\min\{\frac{j}{\alpha},1\}\}} \left\{ V_t^N\left(\frac{j}{\alpha} - z\right) +  (1-z)d_i \right\}\right]  \right. \\
      &  \left. + \sum_{i \in N} \left[\mu^*_i \min_{z \in \left\{0, \min\left\{1, C_N - \frac{j}{\alpha}\right\} \right\}}\left\{  V^{N}_t\left(\frac{j}{\alpha} + z \right)\right\}\right] +  H^{N}\left(\frac{j}{\alpha}\right) \right) \\
      =& \alpha \cdot V_{t+1}^N\left(\frac{j}{\alpha}\right).
        \end{aligned}\end{equation*}

\noindent The first equality holds by (\ref{enerlaatstebewijslabel}). The second equality holds by the induction hypothesis and (\ref{hvergelikingalfa}). The third equality holds by introducing a new variable $z = l/\alpha$. The fourth equality holds 
as $\alpha$ can be taken out. The last equality holds by Definition \ref{valuekie}, using $1 - \sum_{i \in N} \lambda_i^* + \mu_i^* = 0$. By the principle of mathematical induction, this completes the proof. $\hfill \square$

\begin{theorem5} \label{bewijsvoormidpointconvex} Let $f: \mathbb{N}_0 \to 
 \mathbb{R}$ with $f(x) + f(x+2) \geq 2 \cdot f(x+1)$ for all $x \in \mathbb{N}_{0}$ and let $a,b \in \mathbb{N}_{0}$ with $a\leq b$, $c,d \in \{0,2\}$ and $a+c \leq b+d$. Then, it holds that 
\begin{equation*} \min_{\substack{x \in \{a,a+1,\ldots,b\} \\ y \in \{a+c,a+c+1\ldots,b+d\}}} \left\{ f(x) + f(y) \right\} \geq 2\min_{z \in \left\{a+\frac{c}{2},a+\frac{c}{2}+1, \ldots, b+\frac{d}{2}\right\}} \left\{ f(z)\right\}.\end{equation*}
\end{theorem5}

\noindent \underline{Proof of Lemma \ref{bewijsvoormidpointconvex}} \bigskip

\noindent Let $x,y \in \mathbb{N}$. Observe that
\begin{equation*} \begin{aligned} &\min_{\substack{x \in \{a,a+1,\ldots, b\} \\ y \in \{a+c,a+c+1,\ldots, b+d\}}} \left\{ f(x) + f(y) \right\} \hspace{30mm} \\
&\geq \min\left[ \left\{f(z) + f(z+1) \vert z = a + \frac{c}{2},\ldots, b+ \frac{d}{2}-1 \right\}\cup \left\{2f(z) \vert z = a + \frac{c}{2},\ldots, b+ \frac{d}{2} \right\}\right] 
\end{aligned} \end{equation*}
\begin{equation*} \begin{aligned}
& \geq \min\left\{ 2f(z) \hspace{2mm} \vert \hspace{2mm} z = a+\frac{c}{2},\ldots,b+\frac{d}{2}\right\} \hspace{80mm} \\
&=2\min_{z \in \left\{a+\frac{c}{2},a+\frac{c}{2}+1, \ldots, b+\frac{d}{2}\right\}} \left\{ f(z)\right\}.
\end{aligned} \end{equation*}


\noindent The first inequality holds by midpoint convexity, i.e., the fact that for any $x \in \{a,a+1,\ldots,b\}$ and $y \in \{a+c,a+c+1,\ldots, b+d\}$ there exists a $z \in \{ a + \frac{c}{2},\ldots, b+ \frac{d}{2} \}$ such that $f(x) + f(y) \geq f(z) + f(z+1)$ or $f(x) + f(y) \geq 2 f(z)$
, which holds true as long as $f$ is convex (see e.g. \citet[Thm. 1.2]{murota2003discrete}). The second inequality holds as $\min\{2z,z+w,2w\} = \min\{2z,2w\}$ for any $z,w \in \mathbb{R}$. The last equality follows by taking (constant) $2$ outside the minimum.  $\hfill \square$

\end{document}